\documentclass[sigconf, authorversion]{acmart}

\usepackage{booktabs} 

\renewcommand\footnotetextcopyrightpermission[1]{} 
\setcopyright{none}

\acmDOI{10.475/123_4}

\acmISBN{123-4567-24-567/08/06}

\acmConference[]{}{}{}

\begin{document}
\title[Map matching when the map is wrong]{Map matching when the map is wrong: Efficient on/off road vehicle tracking and map learning}

\author{James Murphy}
\affiliation{%
  \institution{Lyft, Inc.}
  \streetaddress{185 Berry Street}
  \city{San Francisco}
  \state{California}
  \postcode{94107}
}
\email{jmurphy@lyft.com}

\author{Yuanyuan Pao}
\affiliation{%
  \institution{Lyft, Inc.}
  \streetaddress{185 Berry Street}
  \city{San Francisco}
  \state{California}
  \postcode{94105}
}
\email{ypao@lyft.com}

\author{Albert Yuen}
\affiliation{%
  \institution{Lyft, Inc.}
  \streetaddress{185 Berry Street}
  \city{San Francisco}
  \state{California}
  \postcode{94107}
}
\email{ayuen@lyft.com}

\renewcommand{\shortauthors}{J. Murphy, Y. Pao, A. Yuen}

\begin{abstract}
Given a sequence of possibly sparse and noisy GPS traces and a map of the road network, map matching algorithms can infer the most accurate trajectory on the road network. However, if the road network is wrong (for example due to missing or incorrectly mapped roads, missing parking lots, misdirected turn restrictions or misdirected one-way streets) standard map matching algorithms fail to reconstruct the correct trajectory. 

In this paper, an algorithm to tracking vehicles able to move both on and off the known road network is formulated. It efficiently unifies existing hidden Markov model (HMM) approaches for map matching and standard free-space tracking methods (e.g. Kalman smoothing) in a principled way. The algorithm is a form of interacting multiple model (IMM) filter subject to an additional assumption on the type of model interaction permitted, termed here as semi-interacting multiple model (sIMM) filter. A forward filter (suitable for real-time tracking) and backward MAP sampling step (suitable for MAP trajectory inference and map matching) are described. The framework set out here is agnostic to the specific tracking models used, and makes clear how to replace these components with others of a similar type. In addition to avoiding generating misleading map matching trajectories, this algorithm can be applied to learn map features by detecting unmapped or incorrectly mapped roads and parking lots, incorrectly mapped turn restrictions and road directions. 
\end{abstract}

%
%



\maketitle

\section{Introduction}
\label{sec:introduction}

Map matching is a process by which a sequence of Global Positioning System (GPS) locations from, for example, a vehicle is matched to a path traversed through a known road network, using map information.  It is sometimes thought of as `snapping' a GPS trace to roads in a map.  This has two distinct applications:  improving the accuracy of received locations, and linking GPS traces to the road network itself. In the first case, map matching can improve the accuracy of location data by using knowledge of the road map to augment direct but possibly noisy sensor data, via the assumptions that vehicles move on roads.  In the second case, map matching allows learning about the road network (for example, learning average road speeds) by linking vehicle trajectories directly to components of the road network.

This paper considers the problem of map matching when the map is incorrect or incomplete. This is important because even the best maps will contain errors, omissions or simply become out of date as the world around them changes.  In that case, the advantages of map matching can become disadvantages: forcing vehicle trajectories to follow a path in an incomplete or incorrect road network may make it \textit{less} accurate rather than more.  If incorrect trajectories are generated it may also lead to learning incorrect information about the state of the world.   

In order to mitigate these issues, a method is introduced here that allows existing sample-based map matching algorithms (e.g. HMM-based methods) to be efficiently made robust to such map issues.  It allows tracking of vehicles both on and off known road networks, switching between the two modes as necessary in a single trajectory.  This allows standard on-road vehicle motion to be map matched as usual, but allows off-road trajectory portions to be generated as a fallback that can be used when the road network as described by the map cannot plausibly accommodate the observed vehicle motion.  This allows map matching to be robust to map errors, helping to improve the accuracy of output trajectories, and giving information about places in which the map is wrong; see, for example, Fig.~\ref{fig:one}.  

The paper is structured as follows. Section~\ref{sec:existing-work} surveys existing work on map matching and on- and off-road vehicle tracking. Section~\ref{sec:simm-outline} outlines the details of the method, with Sec.~\ref{sec:simm-forward-filter} describing a forward filter suitable for realtime vehicle tracking, and Sec.~\ref{sec:simm-MAP-trajectory} describing a backward sampling pass for generating map matched trajectories. Section \ref{sec:results} displays some results of our algorithm, compares them to existing methods for map matching in Sec.~\ref{sec:results-vehicle-localization}, shows that it can be used for map learning in Sec.~\ref{sec:results-map-learning} and analyzes its performance. Section~\ref{sec:conclusions} draws conclusions and makes suggestions for future work.

\section{Existing Work}
\label{sec:existing-work}

Multiple approaches to map matching from GPS trace data have been developed.  Early approaches were often geometric or route based~\cite{Brakatsoulas2005, White2000}; the simplest of these simply match GPS points to the nearest on-road point, but more sophisticated variants add increasing consideration of route consistency constraints and observation plausibility.  Over the last decade hidden Markov model (HMM) approaches have grown in popularity~\cite{Goh2012,Jagadeesh2017, Newson2009, Raymond2012, Thiagarajan2009}.  These approaches represent the vehicle's state (e.g. position, speed, heading, etc.), which cannot be directly observed (i.e. it is hidden).  They then define a model of the system governing the evolution of this hidden state in terms of probabilistic emission and transition models.  The emission model describes the probability of making an observation such as a GPS given a particular vehicle state, whereas the transition model gives the probability of transitioning from one hidden state to another in the next time period.  In this way they can enforce both route consistency and observation plausibility.  In general, the state space must be discretized into finite possible hidden states to allow for tractable inference.  The implementation of map matching in the commonly used Open Source Routing Machine (OSRM) routing engine~\cite{OSRM} is based on the HMM algorithm in~\cite{Newson2009}. Such map matching can also be adapted to multiple transit modes, for example \cite{Chen2015} recently proposed a map matching algorithm able to determine the use of different transit modes.

In addition to map matching work from the GIS community, there is a substantial body of relevant work on Bayesian filtering for object tracking and localization~\cite{Arulampalam2002, Sarkka2007}, largely based on sequential Monte Carlo methods such as the particle filter and its derivatives~\cite{Cappe2007, Doucet2001, VanDerMerwe2001}. This tackles a somewhat similar problem, especially when tracking is augmented with road information, albeit with a focus on realtime filtering and localization accuracy as opposed to trajectory inference, which is the focus of map matching. For standard types of motion, such as that of passenger cars in free-space (i.e. not constrained to move only on roads) closed-form filters such as Kalman filters or approximate nonlinear variants such as the extended or unscented Kalman filter (E/UKF)~\cite{Julier1997, Wan2000} can track targets and infer their trajectories accurately and robustly in a computationally efficient way~\cite{Sarkka2008}.  Such closed-form filters cannot, however, easily deal with the multimodality introduced by the branching structure of road networks.  Sample-based approaches such as particle filtering or finite state space HMMs are more easily able to deal with the nonlinearlity and multimodality induced by restricting states to lie on the road network.  

Our work here draws on both these bodies of work to develop a general approach to map matching that is robust to the road map being wrong.  This uses techniques from road-augmented target tracking, which allow targets to be tracked both on- and off-road, taking advantage of road map data when appropriate to improve tracking accuracy~\cite{Cheng2007, Murphy2014, Orguner2009, Ulmke2006}. These methods use a two-model approach, in which an on-road model and an off-road model are run simultaneously.  At any given time the relative probabilities of the models are calculated and the state distribution taken as a weighted mixture of the output from each.  In order to correctly model certain types of behaviour such as crossing from one road to another via an off-road area, the models must be able to `interact' with each other, so that, going forward in time, new states in one model can be spawned from the current state of the other model.  For example, an off-road state should be able to give rise to a nearby on-road state in the next step of the filter to model the possibility of transitions from off- to on-road movement. A sample-based multiple-model approach able to deal with on- and off-road motion is offered by interacting multiple model particle filter (IMM-PF) approaches, for example in \cite{Orguner2009}, which maintain fixed-size populations of samples for each model and use these to calculate model probabilities at each step.  In these approaches, both on- and off-road tracking is performed with sample-based methods, which can make such methods computationally demanding and lose benefits such as robustness of closed-form tracking methods in the off-road case.

\begin{figure}
\includegraphics[width=3.3in]{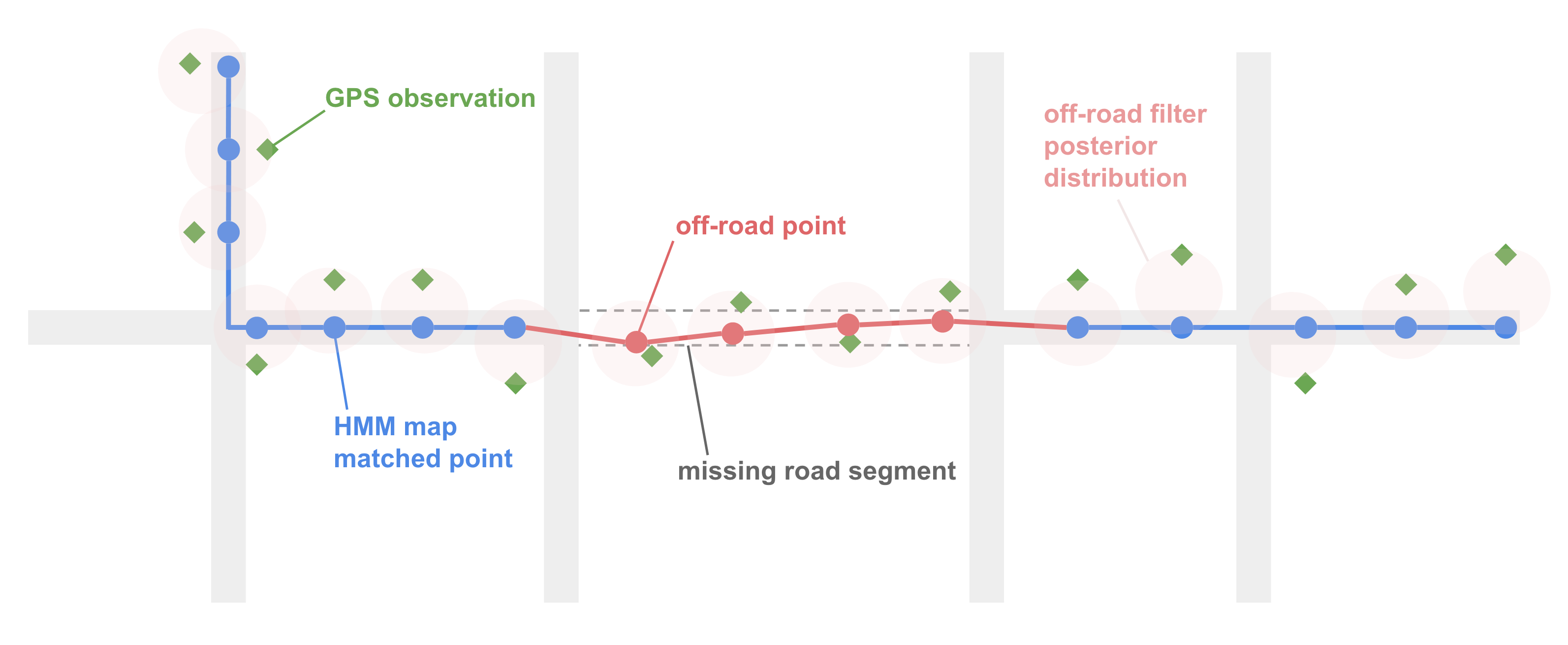}
\caption{Example of idealized output - when the vehicle is moving on correctly mapped roads map matching should be used to produce on-road trajectories.  If a missing segment is encountered a portion of off-road trajectory should be generated; map matching can then continue correctly after the missing segment, and a valid trajectory generated.}
\label{fig:one}
\end{figure}

Here, an approach is outlined which attempts to combine sample-based on-road tracking methods with efficient closed-form off-road tracking.  In particular, the off-road tracking model is intended as a `fallback' model only at times when the on-road model does not provide sufficient flexibility to describe vehicle motion.  For this reason, we would like to spend only a small amount of computation on the fallback off-road filter, and retain maximum flexibility in the on-road map matching model. Figure \ref{fig:one} illustrates the potential benefits of such a system in the case of missing roads. 

The approach developed here is termed the \textit{semi-Interacting Multiple Model} (sIMM) approach.  It is based around the idea that it must be possible to spawn new on-road states from off-road states (see figure \ref{fig:one}), but that it is sufficient for the off-road filter to run independently (hence `semi'-interacting), because it tracks in free space.  This allows a closed-form off-road filter to interact with a sample based method without hypothesis explosion~\cite{Blom1988} and with minimal additional computational cost compared to running the two filters separately.  The assumptions inherent in this framework and their likely effects are made clear in what follows, and are shown to be not especially onerous in this application.  The sIMM framework is presented here in a general way, allowing any sample-based on-road map matching algorithm to be used.  Special attention is given to the finite state space HMM case because of its importance in map matching applications, but alternative models such as particle filters could easily be substituted in this framework.

The sIMM filter can be seen as a special case of the IMM-PF filter in which the off-road tracking filter has a single marginalized sample~\cite{Cappe2007}, and in which interaction is approximated by restricted interaction functions. Taking this view shows how to extend the sIMM model given here to the fully interacting case.

In addition to vehicle tracking and map matching, the sIMM filter can be leveraged to detect errors in the road network. Most past work on map learning has been focused on learning the road network from GPS traces --- e.g., using kernel density estimators~\cite{Biagioni.2012.trr, Biagioni.2012.sigspatial}, iterative methods~\cite{He.2018.sigspatial}, clustering algorithms~\cite{Chen.2016.kdd}, graph spanners~\cite{Stanojevic.2017.arvix} or neural network~\cite{Goyal.2019.kdd}. When GPS traces are not available, with the recent progress in computer vision using deep learning, image segmentation using satellite and aerial imagery has become a popular approach \cite{Bastani.2018.cvpr, Mattyus.2017.iccv}.

\section{Semi-Interacting Multiple Model Filter}
\label{sec:simm}

\subsection{Outline of Approach}
\label{sec:simm-outline}

The key idea behind the sIMM filter is to run a stable, closed-form Kalman filter (or extended or unscented variant, hence denoted (E/U)KF) as a general free-space (off-road) tracking model for the vehicle being tracked.  This filter evolves separately and uninfluenced by an HMM (or other sample-based) on-road tracking component.  Crucially, however, the (E/U)KF off-road tracker \textit{is} allowed to give rise to ancestors of next-generation on-road samples.

This maintains the cheap, robust (E/U)KF off-road filter as an independent component that can be run separately, but allows the on-road tracker to use this to create off-road portions of motion in a principled way.  The (E/U)KF filter does not really need to interact with the on-road filter because, being a free-space tracker, it will not diverge from true vehicle track (at least when observations are somewhat reliable) in the same way that the on-road filter can (e.g. when a link road does not exist, necessitating a large detour).  This is at the cost of some accuracy, since it cannot benefit from map information.  However since this component is considered a fallback, some loss of accuracy is tolerable.

The two models (E/U)KF and HMM each track the target, but it is assumed that the target may switch from motion best described by one model to the other at any time.  The vehicle's overall posterior is taken to be a mixture model over the two tracking models with a component weight assigned to each model.  The interaction between the models comes because transitions are allowed between the models, but here that interaction is `semi-' because the on-road filter is not allowed to influence the evolution of the off-road (E/U)KF filter, but interaction \textit{is} allowed in the opposite direction.

The following section develop the forward semi-interacting multiple model filter (sIMM) and the subsequent section shows how the output of this filter  can be used to generate \textit{maximum a posteriori} (MAP) vehicle trajectories as in map matching applications.

\subsection{Forward sIMM Filter}
\label{sec:simm-forward-filter}

Let $M_t$ be the mode of motion at time $t$, where $M_t \in \{r, g\}$ with $r$ representing the on-road model and $g$ representing the off-road (general) model.  Given a series of observations $y_{1:t}$, the overall model is that the state posterior can be represented as a weighted mixture of the posterior conditioned on each of these models, i.e. that
\begin{align}  			      
p(X_t | y_{1:t}) &= \int p(X_t, M_t | y_{1:t}) dM_t \nonumber\\
&= \int p(X_t | M_t, y_{1:t})p(M_t | y_{1:t}) dM_t\nonumber\\
&= \sum_{m \in \{r, g\}} \mu_t^m p_m(X_t^{m} | y_{1:t})
\end{align}
where $p_m(X_t^{m} | y_{1:t}) = p(X_t | M_t=m, y_{1:t})$ is the \textit{mode-conditioned} posterior of the vehicle state at time-period $t$, i.e. immediately after the $t^\text{th}$ observation, and $\mu_t^m = p(M_t | y_{1:t})$ is the mode weight of mode $m$ at that time.

The mixture weights $\mu_t^m$ can be found as

\begin{align}
\mu_t^m &= \sum_{n\in\{r,g\}}\int p(M_t=m, M_{t-1}=n, X_{t-1}^n, X_t^m | y_{1:t})dX_{t-1}^n dX_t^m \nonumber \\
&\propto \sum_{n\in\{r,g\}} p(M_t=m | M_{t-1}=n) \bigg[ \int p(X^m_t|X^n_{t-1}, M_t=m, M_{t-1}=n)\nonumber\\&\qquad p_m(y_t|X_t^m, M_t=m)p(M_{t-1}=n, X_{t-1}^n | y_{1:t-1}) dX_{t-1}^n dX_{t}^m\bigg]
\end{align}
where the $X^m_t$ represents the state vector of model $m$ at time-period $t$; for example $X_t^r \in R$ and $X_t^g\in G$, where $R$ and $G$ are the state spaces of the on- and off-road tracking models, respectively.

Using a two state Markov chain as a transition model for the motion type, so that the prior probability of going from motion of type $n$ to type $m$ is $\pi_{nm}$, this can be written as 

\begin{align}
\mu_t^m &\propto \sum_{n\in\{r,g\}} \mu_{t-1}^n\pi_{nm} \int p_m(y_t|X_t^m)p_{nm}(X_t^m|X_{t-1}^n)\nonumber\\
&\qquad\qquad\qquad\qquad\qquad\times p_n(X_{t-1}^n | y_{1:t-1}) dX_{t-1}^ndX_t^m
\end{align}
where the notation $p_m$ indicates conditioning on model $m$ at the appropriate time, for example conditioning on $M_t=m$; $p_{nm}$ is used to indicate conditioning on $M_t=m$ and $M_{t-1}=n$.  Thus, the terms in the integral are, in turn: the observation density at time $t$ for model $m$; the state transition density from model $n$ at time $t-1$ to model $m$ at time $t$ (further details are given below); and the state posterior filtering density conditioned on model $n$ at time $t-1$.

Note that the integral term here is the observation likelihood conditional on the motion model in the previous and current time period
\begin{align}  
I_{nm} &= \int p_m(y_t|X_t^m)p_{nm}(X_t^m|X_{t-1}^n)p_n(X_{t-1}^n | y_{1:t-1}) dX_{t-1}^ndX_t^m \nonumber \\
&= p(y_t | y_{1:t-1}, M_t=m, M_{t-1}=n)	         
\end{align}                                    
                                     
The model-conditioned state posterior density (final term in the integral above) is given by the posterior filtering density for the specified motion model.  In the case of the general off-road motion model this will be a Gaussian distribution given by the (E/U)KF, whereas in the case of the on-road model, this is a weighted sample-based distribution given by the HMM.

Evaluating the integral $I_{nm}$ for each combination of $n$ and $m$ is the key challenge in formulating the sIMM, and the rest of this section will outline a series of approximations that can be used in order to do so.  An attempt will be made to make clear the approximations being made and the possible consequences of those approximations.\newline 

\noindent\textbf{Assumption 1:}  Semi-interaction

The eponymous `semi'-interacting approximation amounts to approximating the state transition function when moving from on- to off-road models as
\begin{equation}
p_{rg}(X_t^g | X_{t-1}^r) \approx p_{gg}(X_t^g | X_{t-1}^g)
\end{equation}
It is also assumed that in this case $p_r(X^r_{t-1} | y_{1:t-1}) \approx p_g(X^g_{t-1} | y_{1:t-1})$.  In particular, assuming that the right-hand side is a good approximation of the left in both these approximations, i.e. that the state posterior distribution under the (E/U)KF is similar to that starting from the road state in the previous stage.  

The assumption above will be good whenever the off-road state estimate is already close to the road point under consideration, and will be worse when it is far away.  This is ideal (for an approximation) because when the off-road estimate is far away from the road, the on-road hypothesis should be weak and the $r \to g$ transition will be of minimal importance.  Because the off-road tracking will, in general, fit the observation better than one starting from the on-road point (because its position is not constrained) the probability of the $r \to g$ transition will generally be slightly overestimated, meaning that $r \to g$ transitions will be slightly more likely than they would ideally be.  In the near-to-road case, this effect will be smaller, as the approximation will be better, and in the far-from-road case this effect will be larger, but in this case the probability of an $r \to g$ transition will anyway be small, and therefore the effect of the approximation should be small overall.\newline

\noindent\textbf{Assumption 2:}  Off- to on-road transitions are approximated by assuming the vehicle was at its nearest on-road state at the previous time.

The sIMM here uses the same approximation in the off-road to on-road transition as that used in \cite{Orguner2009}:
\begin{equation}
p_{gr}(X_t^r|X_{t-1}^g) \approx p_{rr}(X_t^r | g2r(X_{t-1}^g))
\end{equation}
where $g2r: G\to R$ is a mapping function from general to road states.  This approximates the off- to on-road state transition density as the road transition density assuming the preceding off-road point was already at $g2r(X_{t-1}^g)$.  This will overestimate the transition density when the off-road point is far from the road, because no account is taken of the distance of the off-road point to the corresponding on-road projection, making this transition density the same for all points in $G$ projecting to $g2r(X_{t-1}^g)$.  This will inflate the $g \to r$ state transition density in the case of off-road points very distant from the road.  In these cases, once again, the overall $g \to r$ transition probability should be small and so the effect of this approximation will be reduced.  An alternative would be to use the full free-space transition model, approximating the transition density as 
\begin{equation}
p_{gr}(X_t^r|X_{t-1}^g) \approx p_{gg}(r2g(X_t^r) | X_{t-1}^g)
\end{equation}
Something closer to an idealized model would be given by
\begin{equation}
p_{gr}(X_t^r|X_{t-1}^g) \approx \int p_{rr}(X_t^r | g2r(X^*_{t^*}))p_{gg}(X^*_{t^*}|X^g_{t-1})dX_{t^*}^*dt^*
\end{equation}
where $X^*_{t^*} \in R$ ranges over the set of road points and $t^*$ is a timestamp between the times of observations $t-1$ and $t$.  However, for practical motion models this will be intractable in reasonable time, so the approximation in equation 6 is used instead.\newline

There are four previous-to-current stage model combinations for which we need to evaluate the model-conditioned likelihood $I_{nm}$ in (4).  Assumptions 1 and 2 define the conditional state transition function in all cases:
\begin{align}
p_{nm}(X_t^m | X_{t-1}^n) = \begin{cases}
p_{gg}(X_t^g | X_{t-1}^g) & m=g, n\in \{r,g\}\\
p_{rr}(X_t^r | X_{t-1}^r) & m=n=r\\
p_{rr}(X_t^r | g2r(X_{t-1}^g)) & m=r, n=g
\end{cases}
\end{align}

The integral $I_{nm}$ can be calculated in each of these cases as follows:

\noindent\textbf{Case $m=g, n \in \{g,r\}$ ($g\to g$ and $r\to g$ transitions):} the integral $I_{nm}$ corresponds to the prediction error decomposition from the (E/U)KF.  This can be calculated in closed form to give
\begin{align}
I_{rg} \approx I_{gg} &= \int p_g(y_t | X_t^g)p_g(X_t^g | y_{1:t-1}) dX_t^g  \nonumber\\
&= \mathcal{N}\left(y_t ; \mu_{y_t}, \Sigma_{y_t}\right)
\end{align}
with, for the extended Kalman filter off-road tracking model,
\begin{align}
\mu_{y_t} &= h_t(\hat \mu_{t|t-1}) \nonumber \\
\Sigma_{y_t} &= H_t \hat\Sigma_{t|t-1}H_t^T + R_t \nonumber \\
H_t &= \frac{\partial h_t}{\partial X}|_{\hat \mu_{t|t-1}}\nonumber
\end{align}
where $h_t(X)$ is the (nonlinear) observation function. $\hat \mu_{t|t-1}$ and $\hat \Sigma_{t|t-1}$ are the predictive state mean and covariance, respectively, calculated in the predict step of the extended Kalman filter. $R_t$ is the covariance matrix of the observation noise.
For the standard Kalman filter, $h_t(X)$ is replaced with $H_tX$, where $H_t$ is the observation matrix, which also replaces the Jacobian in the covariance update.

\noindent\textbf{Case $m=n=r$ ($r \to r$ transition):}  here, the previous-state posterior is given from the sample-based on-road tracker by a weighted sample approximation, i.e. 
\begin{equation}
p_{rr}(X_{t-1}^r | y_{1:t-1}) \approx \sum_j w_{t-1}^{r,j}\delta_{\{ X_{t-1}^{r,j} \}}
\end{equation}
where $\delta_{\{x\}}$ represents a unit point probability mass located at $x$. $X_{t-1}^{r,j}$ is the location of the $j^\text{th}$ sample in the sample-based approximation at observation time $t-1$ and $w_{t-1}^{r,j}$ its associated weight.

Similarly, the current-stage predictive distribution is also given by weighted samples (e.g. in the case of the ad-hoc HMM from \cite{Newson2009} the placement of these is induced by the observation $y_t$).  In this case, both the $X_{t-1}$ and $X_t$ integrals are approximated by finite sums to give
\begin{align}
I_{rr} \approx \sum_i p_r(y_t | X_t^{r,i}) \sum_j p_{rr}(X_t^{r,i} | X_{t-1}^{r,j})w_{t-1}^{r,j} = \sum_i v_t^{r,i}
\end{align}
Note here that the $v_t$ are unnormalized versions of the forward filter weights in the HMM filter, i.e. that
\begin{equation}
w_t^{r,i} \propto v_t^{r,i} \qquad s.t. \qquad \sum_i w_t^{r,i} = 1
\end{equation}

\noindent\textbf{Case $m=r, n=g$ ($g \to r$ transition):}  this is the trickiest case because the previous-stage posterior filter distribution is continuous (Gaussian) but the successor road state is discrete and the mapping $g2r(.)$ is highly nonlinear due to the shape of the road network.  The required integral is given by
\begin{align}
I_{gr} = \sum_i p_r(y_t | X_t^{r,i}) \int p_{rr}(X_t^{r,i} | g2r(X_{t-1}^g))p_g(X_{t-1}^g|y_{1:t-1})dX_{t-1}^g
\end{align}

The inner integral is generally intractable and must be approximated.  This corresponds to the predictive distribution of the road state given that the previous state was off-road and a $g \to r$ transition occurred).  One option is to use a sample-based approximation, approximating the previous filtering posterior with a (possibly weighted) collection of samples.  For example, Monte Carlo methods such as importance sampling, or a deterministic approximation such as a Sigma-Point (unscented) approximation could be used \cite{Julier1997}.  An even easier (but clearly worse!) approximation is to use a single point approximation to the previous posterior filtering state distribution.  Choosing the mean of the previous-state Gaussian filtering posterior distribution $\mu_{t-1|t-1}^g$ gives
\begin{equation}
I_{gr} \approx \sum_i p_r(y_t | X_t^{r,i}) p_{rr}(X_t^{r,i} | g2r(\mu_{t-1|t-1}^g))
\end{equation}      
A better approximation is of course preferable from an accuracy perspective but more computationally expensive.  The (very) simple approximation here uses a single sample at the MAP point of the distribution and so will overestimate this particular likelihood and make $g \to r$ transition more likely than they would otherwise be.  In most cases, this approximation will likely be sufficient, but if $g \to r$ transitions happen too easily, this approximation could be revisited.\newline 

\noindent\textbf{Assumption 3:}  Approximation of previous-state filtering posterior as a single point at the MAP point when calculating $g \to r$ likelihood.\newline

For finite state HMM-based on-road tracking models such as \cite{Newson2009} (that is, models that integrate (sum) over all previous-stage samples, rather than \textit{sampling} the ancestor of each sample as in particle filters) the sample weights must be adjusted at each stage to account for the presence of the off-road model, i.e. to account for the fact that the off-road model can be the ancestor of any on-road sample, and thus contributes to its weight.  The filter-weights for the sample-based on-road model are given as follows:
\begin{align}
w_t^{r,i} \propto \mu_{t-1}^r \pi_{rr} u_t^{r,i} + \mu_{t-1}^g \pi_{gr} p_r(y_t | X_t^{r,i})) p_rr(X_t^{r,i} | g2r(\mu_{t-1|t-1}^g))
\label{eq:filter-weight-adjust}
\end{align}
where $u_t^{r,i}$ is the sample weight calculated in the standard on-road filter without accounting for off-road ancestors, i.e.
\begin{align}
u_t^{r,i} = \sum_j p_r(y_t | X_t^{r,i}) p_{rr}(X_t^{r,i} | X_{t-1}^{r,j}) w_{t-1}^{r,j}
\end{align}
where $j$ ranges over the previous-stage HMM samples.

For particle filter models, this adjustment is not necessary, since the ancestor of each current-stage road sample is \textit{sampled} from the previous stage distribution; this reflects the different approximation made by particle filters compared to finite state HMMs.

\subsubsection{Algorithm}
The overall forward sIMM filtering algorithm, assuming a finite state space-based HMM on-road tracking model, is given as follows (for particle filter-based on-road tracking, this should be appropriately adjusted to reflect the more direct computation of the model weights, see e.g. \cite{Orguner2009}):
\begin{enumerate}
\item Initialization
\begin{enumerate}
\item Initialize the on-road filter as a collection of on-road samples
\item Initialize the off-road filter
\end{enumerate}
\item For each observation $y_t$ with $t = 1,..., T$:
\begin{enumerate}
\item Update the off-road filter to find the posterior $p_g(X_t^g | y_{1:t}) $
\item Calculate the off-road ($g \to g$) observation likelihood 
\begin{displaymath}
I_{gg} = p_{gg}(y_t | y_{1:t-1})
\end{displaymath}
using equation 10; this also serves as an approximation to the likelihood $I_{rg}$ in the ($r \to g$) case.
\item Update the on-road filter to find the posterior $p_r(X_t^r| y_{1:t})$
\item Calculate the on-road ($r \to r$) observation likelihood using equation 12:
\begin{displaymath}
I_{rr} = p_{rr}(y_t | y_{1:t-1})
\end{displaymath}
\item Calculate the observation likelihood in the ($g \to r$) case as $c$ using the sum in equation 15:
\begin{itemize}
\item Find the corresponding on-road point of the previous-stage filter posterior mode as $g2r(\mu_{t-1|t-1})$
\item Summing over all samples $i$ in the HMM state approximation at stage $t$:
\begin{displaymath}
I_{gr} = \sum_i p_r(y_t | X_t^{r,i}) p_{rr}(X_t^{r,i} | g2r(\mu_{t-1|t-1}^g))
\end{displaymath}
\end{itemize}

\item If using a finite state space HMM on-road model, adjust the filter's sample weights to those given in equation \eqref{eq:filter-weight-adjust}

\item Calculate the model probability as proportional to:
\begin{align}
m_t^r &= \mu^r_{t-1} \pi_{rr} I_{rr}  + \mu^g_{t-1} \pi_{gr} I_{gr} \nonumber\\
m_t^g &= \mu^r_{t-1} \pi_{rg} I_{rg}  + \mu^g_{t-1} \pi_{gg} I_{gg} \nonumber
\end{align}
\item Normalize for forward filter model probability at stage $t$:
\begin{align}
\mu_t^r &= m_t^r / (m_t^r + m_t^g)\nonumber\\
\mu_t^g &= 1 - \mu_t^r \nonumber
\end{align}
\item The overall filter posterior at stage t is given by:
\begin{equation}
p(X_t | y_{1:t}) = \mu_t^g  p_g(X_t^g | y_{1:t}) + \mu_t^r  p_r(X_t^r | y_{1:t})\nonumber
\end{equation}
\end{enumerate}
\end{enumerate}

\subsection{MAP State Sequence Estimation}
\label{sec:simm-MAP-trajectory}

The filter above allows tracking of targets moving on- and off- the known map.  However, map matching requires inferring the most likely trajectory taken by the target, i.e. the \textit{maximum a posteriori} (MAP) location sequence.  The MAP state sequence $M^*_{1:T}$, $X^*_{1:T}$ is chosen so as to maximize the all-data posterior probability distribution $p(M_{1:T}, X_{1:T} | y_{1:T})$.  By decomposing this distribution in a `cascade', it can be seen how to sample it sequentially to obtain the MAP sequence:
\begin{equation}
p(M_{1:T}, X_{1:T} | y_{1:T}) = p(M_T, X_T | y_{1:T}) \prod_{t=1}^{T-1} p(M_t, X_t | M_{t+1}, X_{t+1}, y_{1:t})
\end{equation}

Here, the distribution $p(M_t, X_t | M_{t+1}, X_{t+1}, y_{1:t)}$ depends only on subsequent state variables (and preceding observations).  The final distribution $p(M_T, X_T | y_{1:T})$ is the final stage filtering posterior distribution, as calculated above.  We can thus sample the MAP point from this by choosing the mode of the posterior filtering distribution at the final stage (discussed for the sIMM model below).  This then allows recursive backward sampling of the MAP state and model going from stage $t = T-1,..., 1$ choosing the mode model and state denoted $M_t^*$ and $X_t^*$ of the distribution $p(M_t, X_t | M^*_{t+1}, X^*_{t+1}, y_{1:t})$ at each stage $t$ and using this in subsequent steps.

In the sIMM model, the meaning of the \textit{mode} of the state and model distribution is not well-defined because the inferred state distribution is a mixture of continuous and discrete probability distributions.  These are not directly comparable (because in the latter case we calculate a probability of being at each of a finite number of points whereas in the former we calculate a probability \textit{density} over the entire state space).   However, it is possible to perform backwards sampling (e.g. for MAP sequence estimation) by first sampling the model indicator and then sampling the state itself from the corresponding model's conditional state distribution.  This corresponds to sampling from $p(M_t, X_t | M^*_{t+1}, X^*_{t+1}, y_{1:t})$ using the cascade decomposition
\begin{align}
p(M_t, X_t | M^*_{t+1}, X^*_{t+1}, y_{1:t}) &= p(X_t | M_t, M^*_{t+1}, X^*_{t+1}, y_{1:t}) \nonumber \\ & \qquad \qquad \times p(M_t | M^*_{t+1}, X^*_{t+1}, y_{1:t})
\end{align}

The marginal distribution of $M_t$ can be found by integrating the above distribution over $X_t$, to give
\begin{align}
p&(M_t | M^*_{t+1}, X^*_{t+1}, y_{1:t}) = \int p(M_t, X_t | M^*_{t+1}, X^*_{t+1}, y_{1:t}) dX_t \nonumber \\
&\propto \int p(X_{t+1}^* | X_t, M_t, M_{t+1}^*)p(M^*_{t+1} | M_t) p(X_t|M_t, y_{1:t})p(M_t | y_{1:t}) dX_t \nonumber \\
&\propto \mu_t^{M_t} \pi(M_t, M_{t+1}^*) \int p(X_{t+1}^* | X_t, M_t, M_{t+1}^*) p(X_t|M_t, y_{1:t}) dX_t
\label{eq:19}
\end{align}
where $\mu_t^{M_t}$ is the model posterior for model $M_t \in \{r, g\}$ calculated in the forward filter, and $\pi(M_t, M_{t+1}^*)$ is mode transition probability from model $M_t$ to $M_{t+1}^*$. We assumed the target motion to be Markovian.
Evaluating the integral above requires the previous filtering posterior for both models (the second term) along with the conditional transition model for the state in each of the four cases $\{r, g\} \to \{r, g\}$.  The same-model transitions ($r \to r$, $g \to g$) are simply those from the corresponding filter.  In the filtering section above, the $r \to g$ and $g \to r$ transition models were both defined.  We could use these, but here we take the slightly unusual step of making different approximations for each of these transitions in the backward direction.

\noindent\textbf{Case $r \to g$:}  This case is the transition from an on-road point at stage $t$ (i.e. $M^*_t = r$), to an off-road point already sampled at $t+1$ (i.e. $M^*_{t+1} = g$).  The forward transition in this case was approximated using the semi-interacting assumption in equation 5.  In backward sampling, if a mapping function $r2g: R \to G$ is available from on-road to off-road states, this earlier approximation can be somewhat atoned for by using the transition density that we would have used in a fully interacting IMM model: 
\begin{equation}
p_{rg}^*(X_{t+1}^{*,g} | X_t^r) \approx p_{gg}(X_{t+1}^{*,g}| r2g(X_t^r))
\end{equation}
This requires that the map $r2g: R \to G$ must encompass all dynamic state components of $G$, which can be tricky and necessitate further approximation when going from on-road states with no or limited dynamic information as in the HMM model in \cite{Newson2009}.

This lets at least some account be taken of the location (and possibly other dynamic state) of the presumed road exit point when sampling the previous stage's road position.  This must be evaluated for each candidate on-road sample at stage $t$, and the MAP point can be chosen by choosing that which maximizes equation \eqref{eq:19} in this case.

\noindent\textbf{Case $g \to r$}:   Here the transition is from an off-road point at stage t (i.e. $M_t = g$), to an already sampled on-road point $X^*_{t-1}$ at $t+1$ (i.e. $M^*_{t+1} = r$).  The sampled on-road point's corresponding off-road point $r2g(X^*_{t+1})$ could thus be as the subsequent state in an off-road model transition.  This approximates the transition from $t$ to $t+1$ as taking place entirely in the off-road model, i.e. making the approximation during the backward pass
\begin{equation}
p_{gr}^*(X_{t+1}^{*,r} | X_t^g) \approx p_{gg}(r2g(X_{t+1}^{*,r})| X_t^g)
\end{equation}
This gives a backward sampling step exactly as that in the $g \to g$ case, replacing $X^{*,g}_{t+1}$ with $r2g(X^{*,r}_{t+1})$.

\begin{figure}
\includegraphics[width=3.2in]{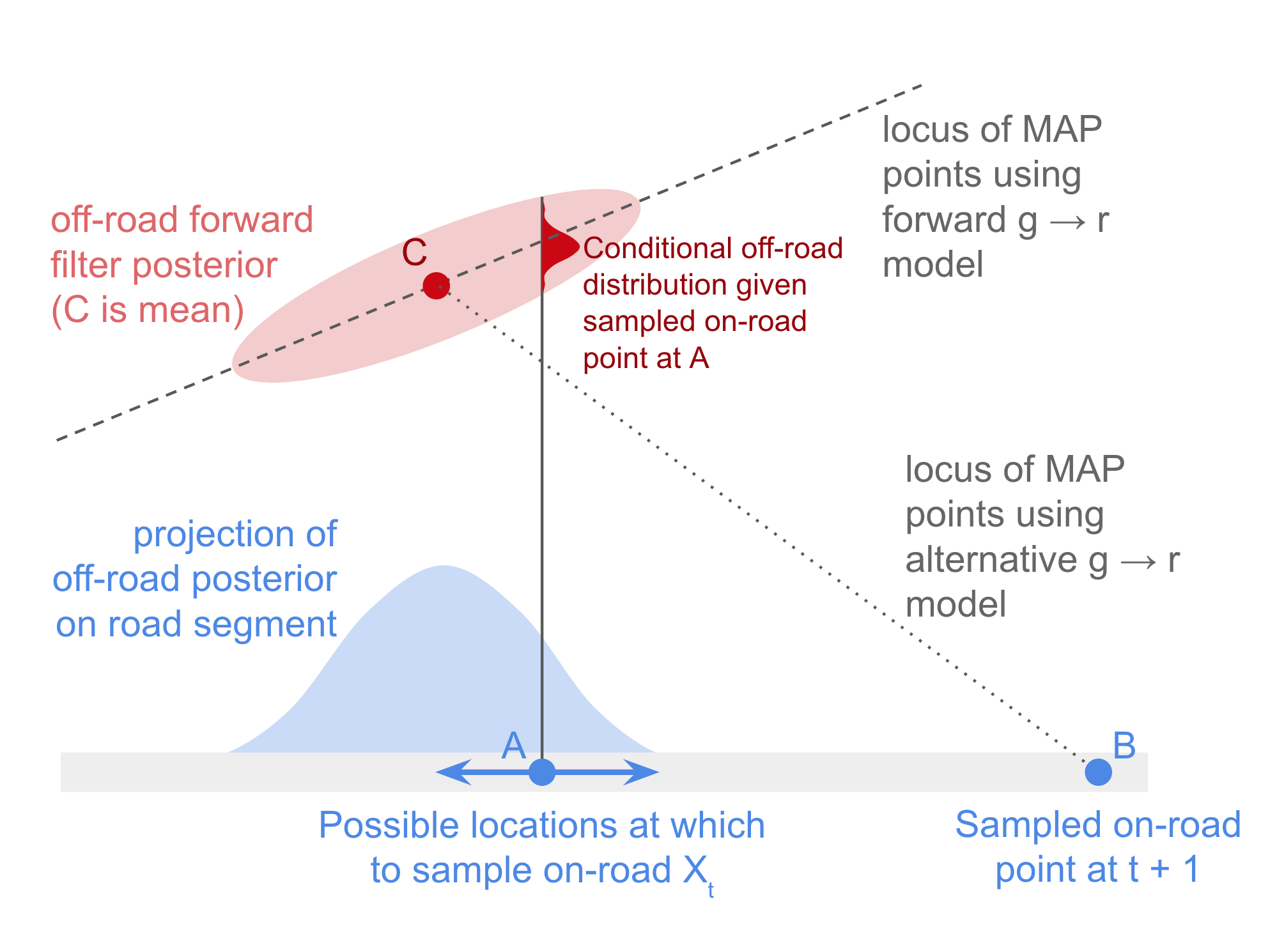}
\caption{Example of how using the $g \to r$ transition model used in the forward direction is not ideal during backward sampling case.  Because it makes use of projection onto the road and the on-road transition model, the forward transition model would result in first sampling an on-road point at $t$, then mapping it to an off-road point, giving a point on the dashed line (assuming nearest road point projection).  However, using the alternative backward $g\to r$ transition model based on the off-road transition model described in section \ref{sec:MAP-trajectory} a point on the line C-B will be generated (in the case of obvious motion models).}
\label{fig:two}
\end{figure}

The backward transition approximation in equation 23 is different from the approximation for $r \to g$ transitions used in the forward direction, and thus adds an additional and non-obvious approximation to the system.  The forward transition model in this case could be used via a sampling approach.  However, this is problematic because of the projection $g2r(.)$ used to map off-road points to on-road points.  Under simple versions of that, which, for example, map off-road points to their nearest on-road point, every point perpendicular to the same on-road point (and not closer to some other road) will map to the same on-road location.  It can be shown that this gives rise to a sometimes poor and counter-intuitive set of possible optima in the backward posterior, that do not account well for the subsequent (already sampled) on-road position; see figure \ref{fig:two}.  \newline

\noindent\textbf{Assumption 4:} In the backward direction, a transition from off-road to on-road is modelled as an off-road transition to the equivalent off-road point.

\noindent\textbf{Effect:}  Distorts the off-road all-data MAP point, especially near $g \to r$ transition.\newline

Given these new assumptions about the transition model in the backwards sampling phase, the integral in equation \eqref{eq:19} can be evaluated in every case, as follows.  Let
\begin{equation}
J_{n\to m} = \int p(X_{t+1}^{*,m} | X_t^n, M_t=n, M_{t+1}^*=m) p(X_t^n|M_t=n, y_{1:t}) dX_t^n. 
\end{equation}
The four cases of this integral can be evaluated as follows:

\noindent\textbf{Case $g \to g$:}
\begin{align}
J_{g\to g} &= \int \mathcal{N}\left(X_{t+1}^{*,g}; f(X_t^g), Q_t\right)\mathcal{N}\left(X_t^g; \mu_{t|t}, \Sigma_{t|t}\right) dX_t \nonumber \\
&= \mathcal{N}\left(X_{t+1}^{*,g}; f(\mu_{t|t}), Q_t + F \Sigma_{t|t}F^T\right)
\end{align}
where $f$ is the state transition function, $F$ is its Jacobian, and $Q_t$ is the state transition covariance from the chosen extended Kalman off-road filter. $\mu_{t|t}$ and $\Sigma_{t|t}$ are the mean and covariance of the posterior filtering distribution from that filter at stage $t$.

\noindent\textbf{Case $r \to g$:}
\begin{align}
J_{r\to g} \approx \sum_i p_{gg}(X_{t+1}^{*,g} | r2g(X_t^{r,i}))w_t^{r,i}
\end{align}

\noindent\textbf{Case $r \to r$:}
\begin{align}
J_{r\to r} = \sum_i p_{rr}(X_{t+1}^{*,r} | X_t^{r,i})w_t^{r,i}
\end{align}

\noindent\textbf{Case $g \to r$:}
\begin{align}
J_{g\to r} &\approx \int \mathcal{N}\left(r2g(X_{t+1}^{*,r}); f(X_t^g), Q_t\right)\mathcal{N}\left(X_t^g; \mu_{t|t}, \Sigma_{t|t}\right) dX_t \nonumber \\
&= \mathcal{N}\left(r2g(X_{t+1}^{*,r}); f(\mu_{t|t}), Q_t + F \Sigma_{t|t}F^T\right)
\end{align}

The MAP model $M^*_t$ can be sampled by evaluating the appropriate versions of the integral (depending on $M^*_{t+1}$) and then choosing the value of $M_t$ that maximizes the expression in equation \eqref{eq:19}.   

Once $M_t^*$ has been sampled the MAP state, $X_t^*$ can be sampled as the mode of the distribution
\begin{align}
p(X_t | X_{t+1}^*, M_t^*, M_{t+1}^*, y_{1:t}) &\propto   p(X_{t+1}^* | X_t, M_t^*, M_{t+1}^*) p(X_t|M_t^*, y_{1:t})
\end{align}
The second term here is the filter posterior distribution and the first is the conditional state transition density.  In order to sample $X^*_t$, each type of $M^*_t$ to $M^*_{t+1}$ transition must be considered.

\noindent\textbf{Case $r \to r$:}  The easiest case to deal with is when $M^*_t = M^*_{t+1} = r$, since this is a standard backward sampling step from the HMM.
\begin{align}
p(X_t^r | X_{t+1}^{*,r}, M_t^*=M_{t+1}^*=r, y_{1:t})  &\propto   p_{rr}(X_{t+1}^{*,r} | X_t^r) p_r(X_t^r| y_{1:t})\nonumber \\
&= \sum_i w_t^{r,i}p_{rr}(X_{t+1}^{*,r}|X_t^{r,i}) \delta_{\{X_t^{r,i}\}}
\end{align}
The MAP $X^*_t$ is then simply the sample $i$ with the maximum backward sampling weight $w_t^{r,i}p_{rr}(X_{t+1}^{*,r}|X_t^{r,i})$.

\noindent\textbf{Case $r \to g$:} The relevant distribution is given (approximately) by
\begin{align}
p(X_t^r | X_{t+1}^{*,r}, M_t^*=r, M_{t+1}^*=g, y_{1:t}) \propto \sum_{i} p_{gg}&(X_{t+1}^{*,g} | r2g(X_t^{r,i})) \nonumber \\
  & \times w_t^{r,i}\delta_{\{X_t^{r,i}\}}
\end{align}
Again, the MAP $X^*_t$ is the sample $i$ with the maximum backward sampling weight. 

\noindent\textbf{Case $g \to g$:}  In the $M^*_t = M^*_{t+1} = g$ case the backward sampling step is that from standard backward sampling in the (E/U)KF and is given as:
\begin{align}
p(X_t^g | X_{t+1}^{*,g}, M_t^*=g, M_{t+1}^*=g, y_{1:t}) &\propto p_{gg}(X_{t+1}^{*,g} | X_t^g) p_g(X_t^g| y_{1:t}) \nonumber \\
&\propto \mathcal{N}\left(X_t^g; \mu_t^*, \Sigma_t^*\right)
\end{align}
with
\begin{align}
\Sigma_t^* &= \left(\Sigma_{t|t}^{-1} + F^TQ_t^{-1}F\right)^{-1}\nonumber \\
\mu^*_t &= \Sigma^*_t\left( \Sigma_{t|t}^{-1} \mu_{t|t} + F^TQ_t^{-1}X_{t+1}^{*,g}\right)\nonumber
\end{align}
The MAP $X^*_t$ is then given by $\mu^*_t$.

\noindent\textbf{Case $g \to r$:}  Using the earlier approximation, this case is essentially the $g \to g$ case above, with $X_{t+1}^{*,g}$ replaced with $r2g(X_{t+1}^{*,r})$, giving
\begin{align}
p(X_t^g | X_{t+1}^{*,r}, M_t^*=g, M_{t+1}^*=r, y_{1:t}) &\propto \mathcal{N}\left(X_t^g; \mu_t^*, \Sigma_t^*\right)
\end{align}
with $\Sigma_t^*$ as above, and 
\begin{align}
\mu^*_t &= \Sigma^*_t\left( \Sigma_{t|t}^{-1} \mu_{t|t} + F^TQ_t^{-1}r2g(X_{t+1}^{*,r})\right)\nonumber
\end{align}
This completes all necessary components to sample the MAP path from the sIMM tracker.

\begin{figure*}
\includegraphics[height=1.7in]{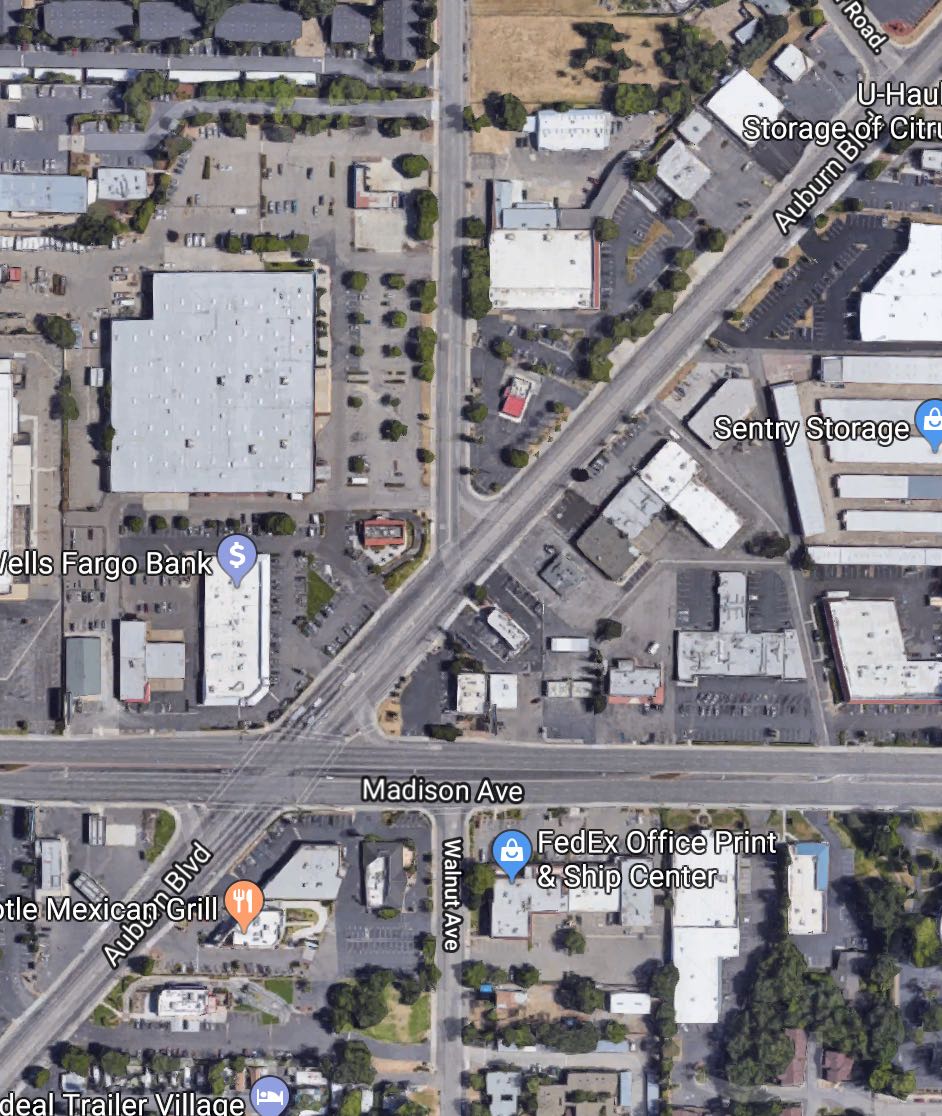}
\quad
\includegraphics[height=1.7in]{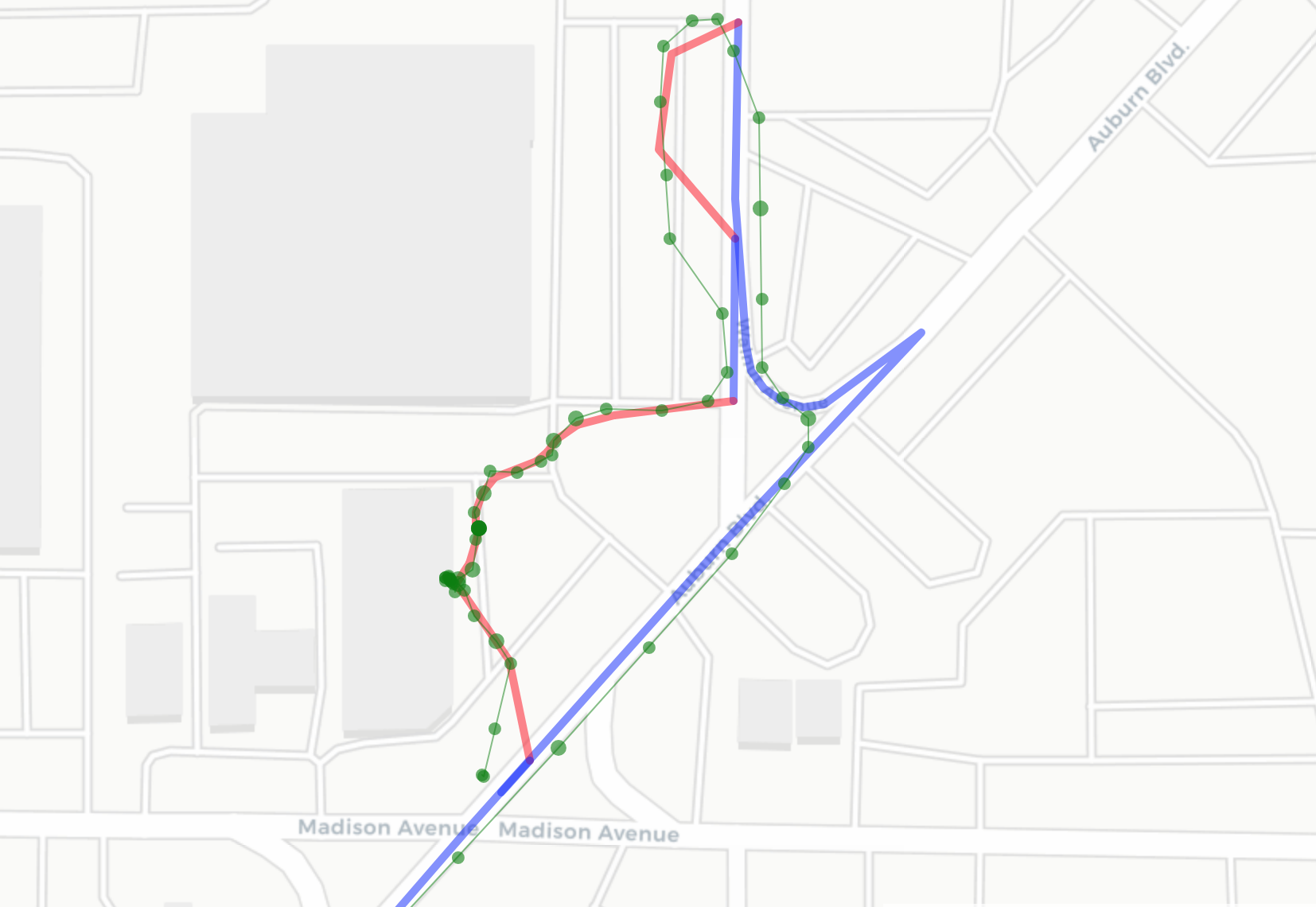}
\quad
\includegraphics[height=1.7in]{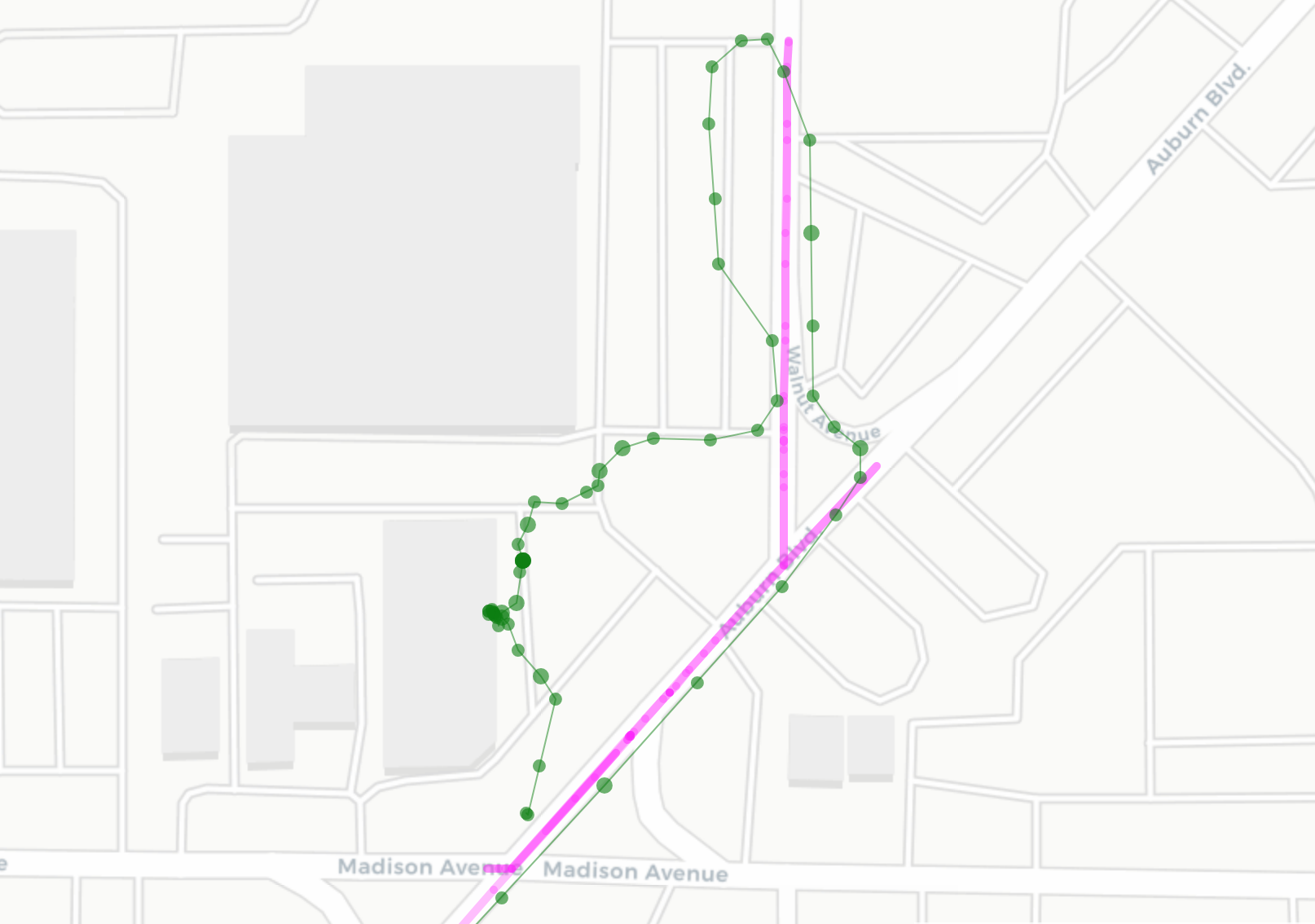}
\caption{Example of on/off-road map matching for a short gps trajectory.  The panel on the left shows an aerial view of the region (taken from Google Maps); the central panel shows the gps trace as green points, on-road map-matched sections as blue lines, and off-road sections of the final trajectory as red lines, computer by the sIMM map matcher. For comparison, the right-hand panel shows the output of standard map matching based on~\cite{Newson2009} with the magenta line showing the map matched trajectory.}
\label{fig:res1}
\end{figure*}

\begin{figure}
\includegraphics[width=3.5in]{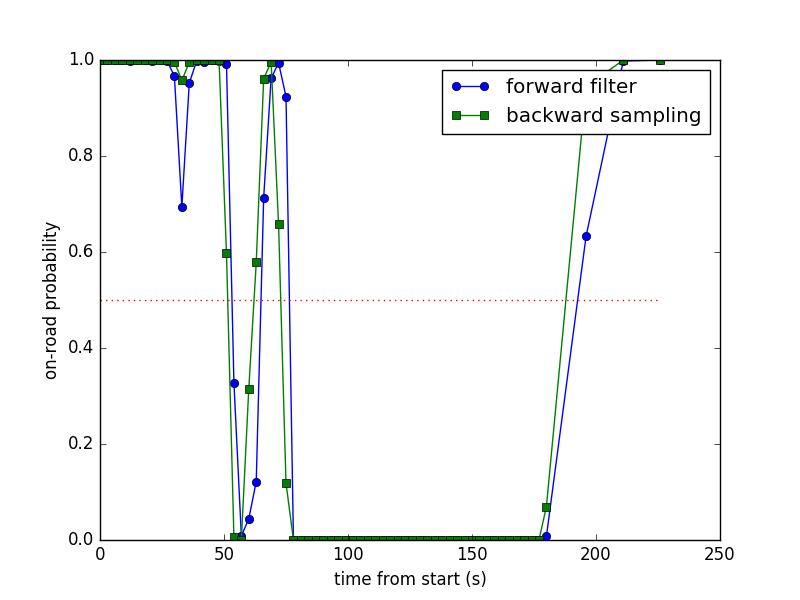}
\caption{On-road probability $\mu_t^r$ generated by the forward sIMM filter (blue circles) and during backward sampling $p(M_t | M^*_{t+1}, X^*_{t+1}, y_{1:t})$ (green squares) for the GPS trace in figure \ref{fig:res1}.}
\label{fig:res1b}
\end{figure}

\subsubsection{Algorithm}
\label{sec:simm-algorithm}
The complete algorithm for backward sampling of the MAP trajectory is given as follows:

\begin{enumerate}
\item Run the forward filtering algorithm above and store:
\begin{itemize}
\item The set of on-road samples and their forward-filter weights $w_t^{r,i}$ at each stage
\item The filter model posterior probabilities at each stage $\mu_t^m$ for $m \in \{r, g\}$
\item The filter posterior distribution for the off-road filter at each stage (fully specified by $\mu_{t|t}$ and $\Sigma_{t|t}$ from the off-road filter) 
\end{itemize}
\item Sample the MAP state at stage $T$ from the final (stage $T$) filter posterior as:
\begin{align}
M_T^* &= \text{argmax}_m \mu_T^m \nonumber \\
X_T^* &= \text{argmax}_{X_T} p_{M_T^*}(X_T | y_{1:T}) \nonumber
\end{align}
\item Iterating backwards for $t = T-1,...,1$:
\begin{enumerate}
\item Evaluate equation \eqref{eq:19} for both values of $M_t$, using the appropriate two cases from equations 25-28
\item Sample $M^*_t$ as that which gives the maximum value in step a. above
\item Given $M^*_t$, use the appropriate sampling strategy from equations 30-33 to sample the MAP state $X^*_t$
\end{enumerate}
\end{enumerate}

\section{Results}
\label{sec:results}

\subsection{Vehicle Localization}
\label{sec:results-vehicle-localization}

\begin{figure*}
\includegraphics[width=2in]{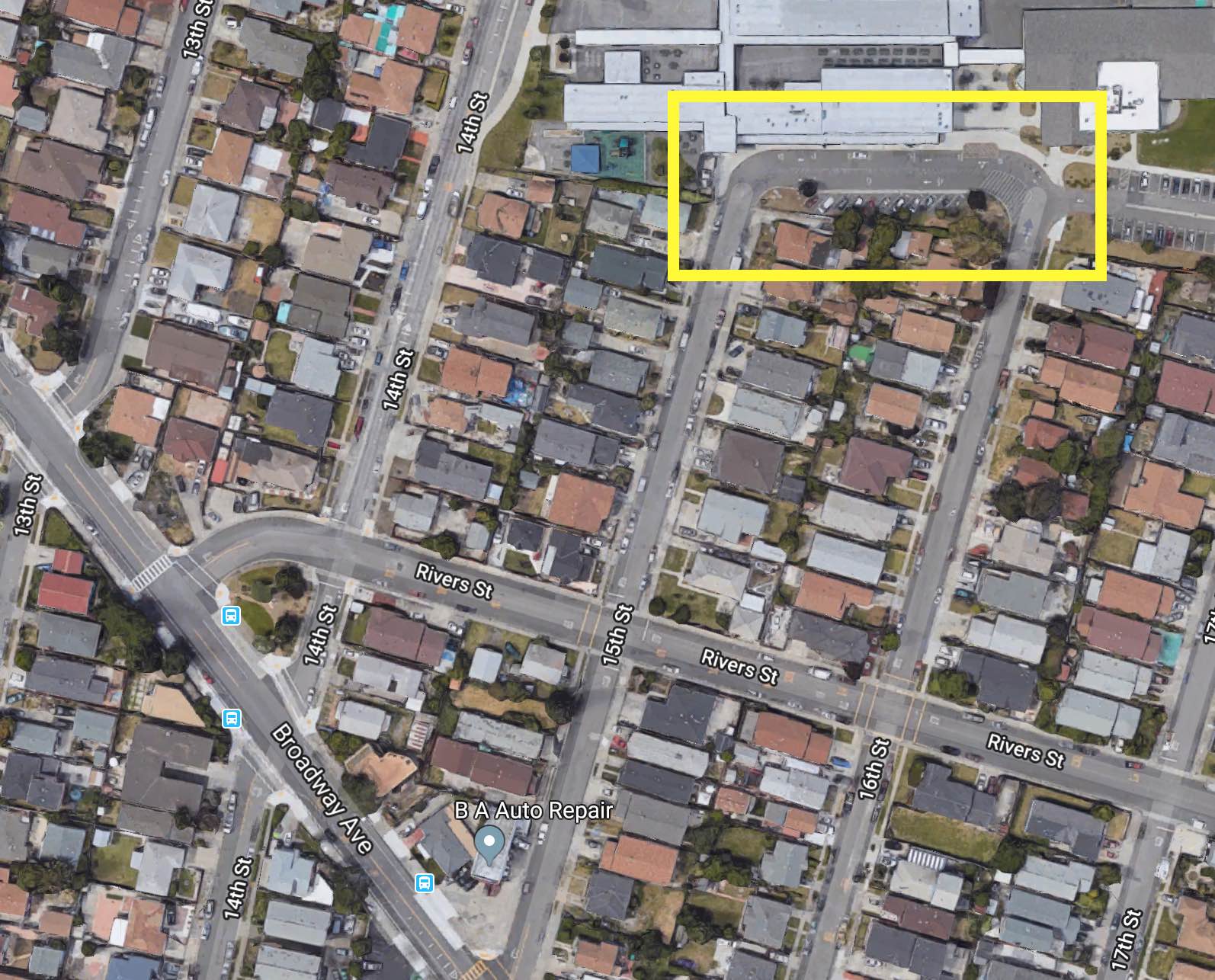}
\qquad
\includegraphics[width=2in]{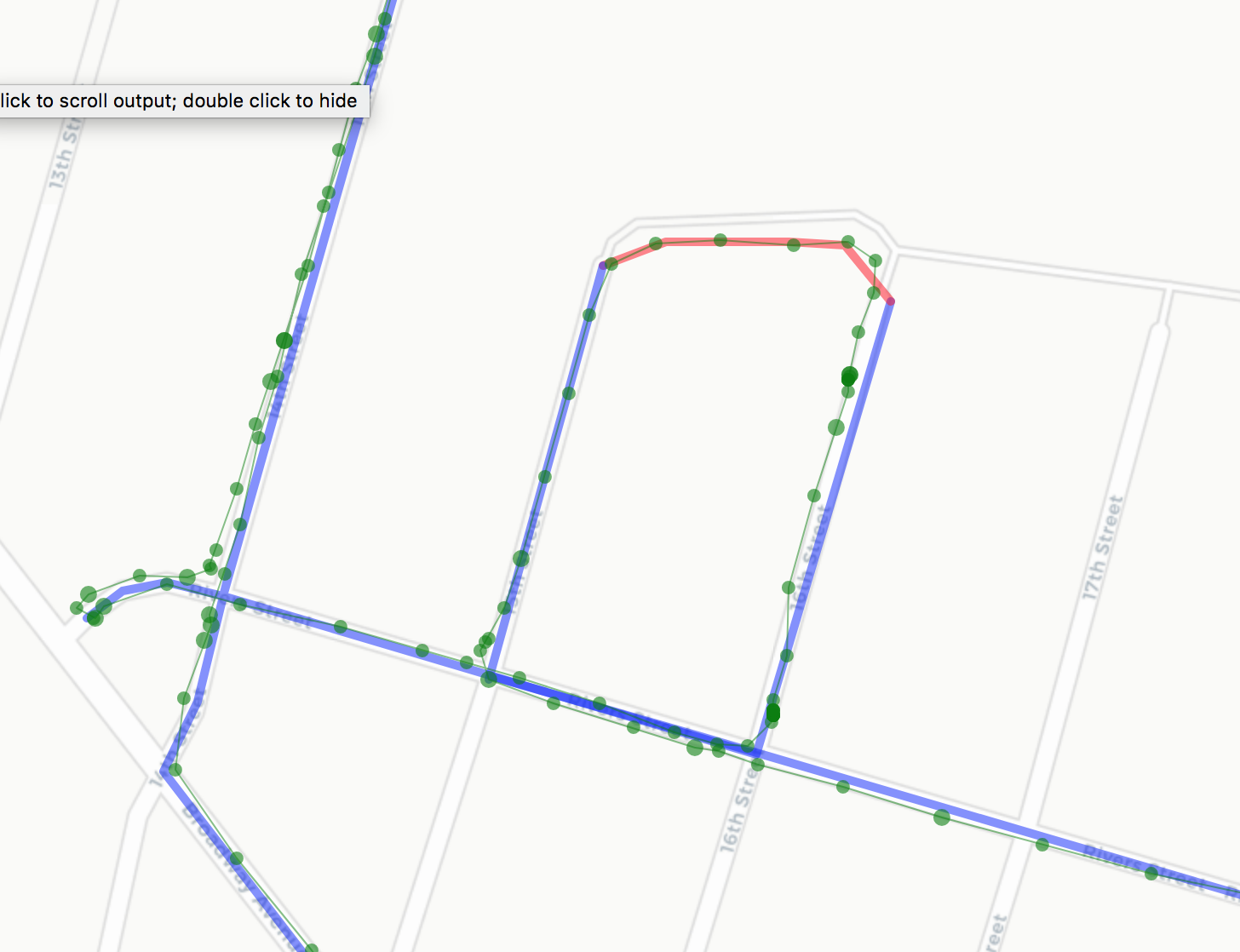}
\qquad
\includegraphics[width=2in]{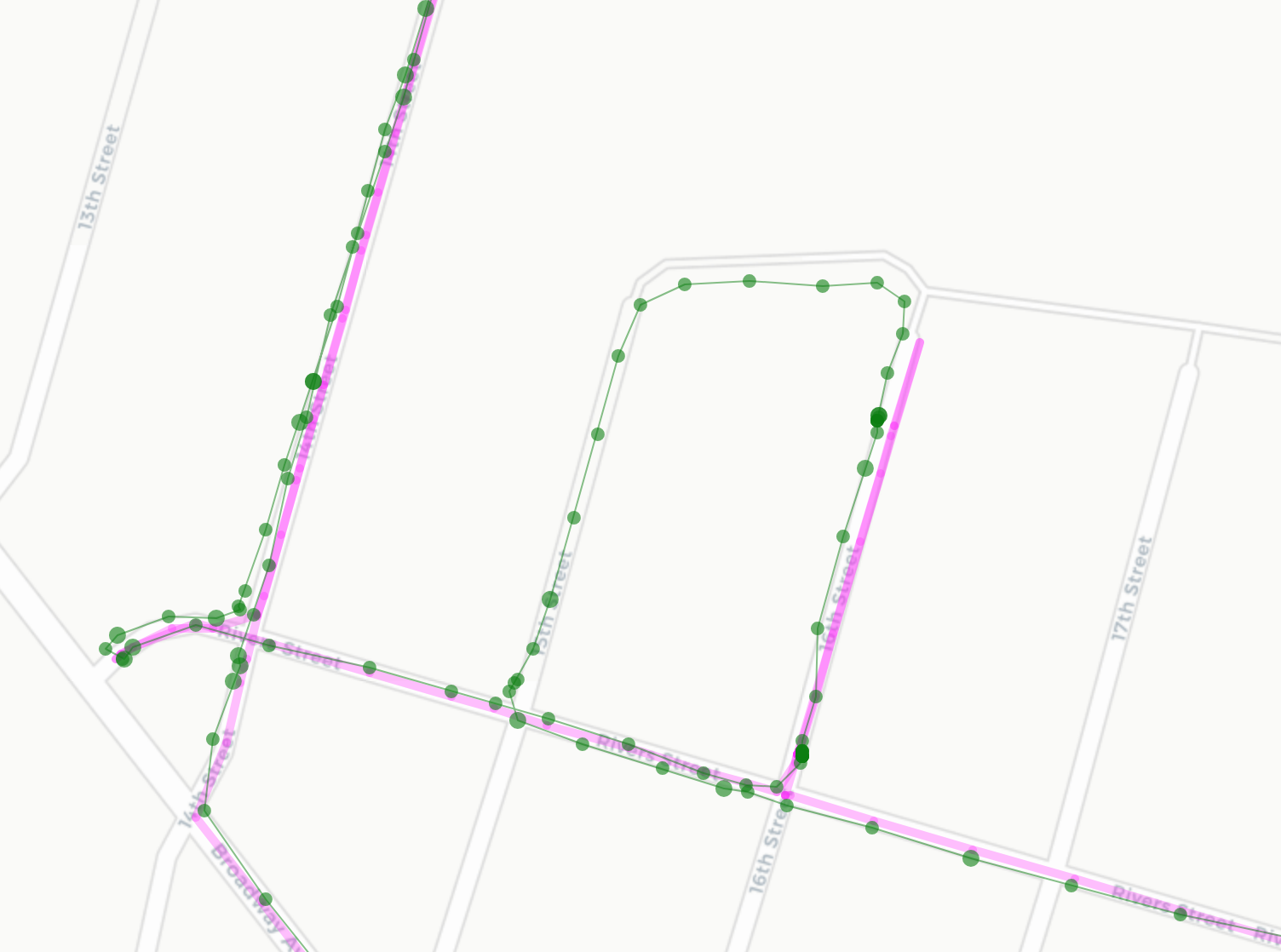}
\caption{Example of on/off-road map matching showing robustness to missing map data.  The road highlighted in the satellite image (from Google Maps) in the left panel is not considered `routable'.  This leads to incorrect map matching using a standard map matching approach (right panel).  On-/off-road map matching (centre panel) correctly traces the route, highlighting the `missing' navigable road.}
\label{fig:res2}
\end{figure*}

\begin{figure}
\includegraphics[height=2in]{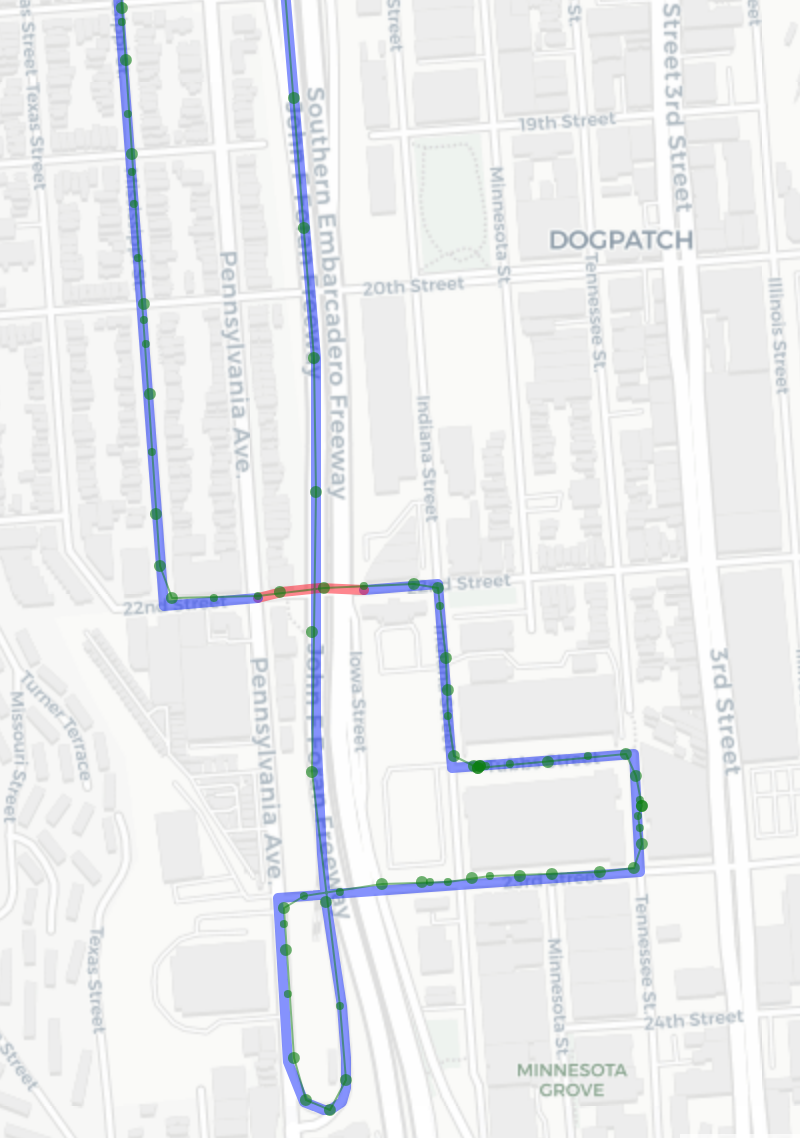}
\qquad
\includegraphics[height=2in]{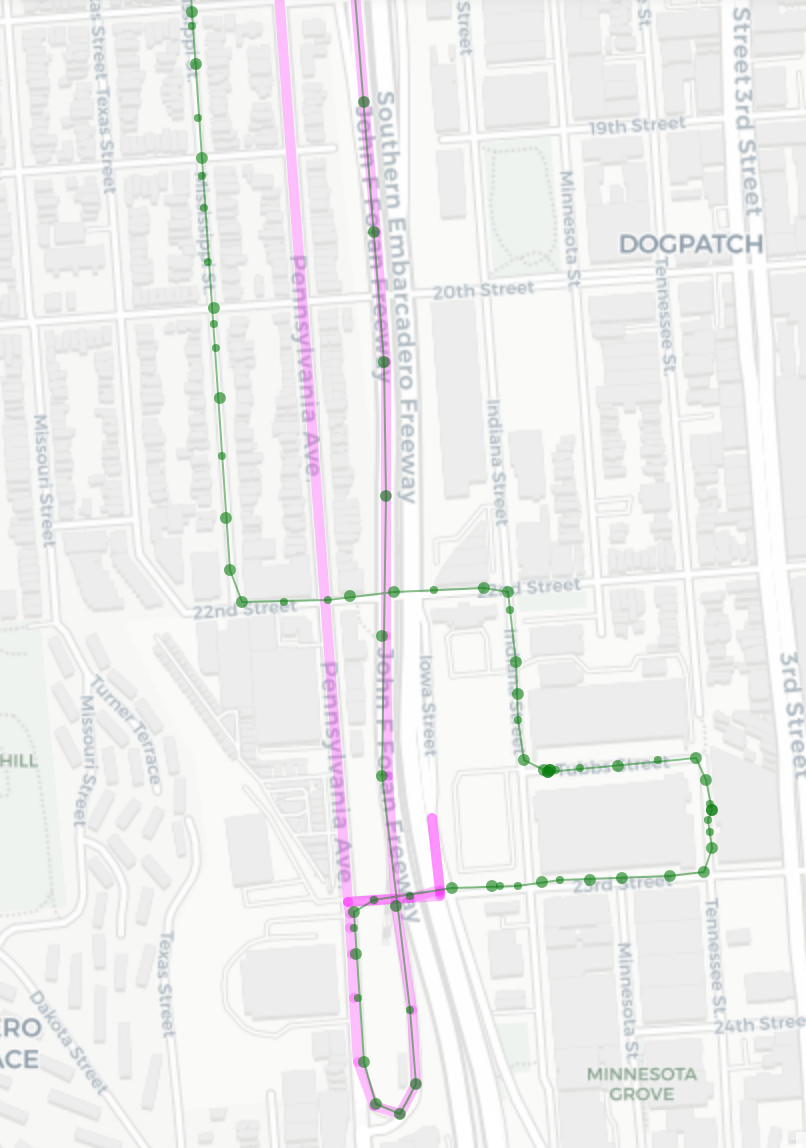}
\caption{On-/off-road map matching with missing/mislabeled road (marked as non-navigable in map).  Left panel shows MAP trajectory output from on-/off-road map matching (on-road portions are blue, off-road are red); right panel shows that from standard road-constrained map matching}
\label{fig:res3}
\end{figure}

Figures~\ref{fig:res1},~\ref{fig:res1b},~\ref{fig:res2} and ~\ref{fig:res3} show the result of running the proposed sIMM filter and MAP trajectory inference on GPS traces for which the road map does not adequately capture the set of drivable places. The map used here is a variant of Open Street Map (OSM) \cite{OpenStreetMap}, dated January 1st, 2018.

Figure~\ref{fig:res1} compares the output of the on-/off-road map matching proposed here to that of a standard (on-road only) map matching approach. The latter is based on a variant of the Open Source Routing Machine (OSRM), which, in turn, is based on the algorithm in \cite{Newson2009}. In this trajectory, based on about four minutes of GPS data sampled every 3 sec, the car, heading northeast from the bottom of the frame, first makes an unusually late turn on to the northerly street, then enters an unmapped parking lot.  The on-/off-road map matcher successfully tracks this motion, outputting a plausible trajectory.  The road-constrained map matcher, in contrast, outputs an implausible trajectory containing several unlikely U-turns.  Figure~\ref{fig:res1b} shows the on-road probability $\mu_t^r$ calculated in the forward sIMM filter. From this, it can be seen that the filter is mostly confident of being in either the on- or off-road state, with uncertainty limited to periods near transition times.  The unusual late turn (around 30s in) shows up as slight possibility of off-road motion in the filter, since it corresponds to a time at which the car's actual motion deviated from that `permitted' on the road network.  Comparing calculated filter mode probabilities to those from the backward MAP trajectory reconstruction pass shows how the backward pass reduces lag in transitions; whereas the forward filter may need to see a couple of observations before it is confident in a mode transition, the backward pass, having the benefit of future information can often make the transition a stage or two earlier, seen in the left-shift of the green curve in Fig.~\ref{fig:res1b}.

Figure~\ref{fig:res2} illustrates how missing roads are handled by the proposed on-/off-road map matcher. In this case the east-west road (highlighted in the first panel) is not marked as drivable in the map, but in fact it is a paved roadway, sufficiently wide for driving.  Because the road-constrained map matcher cannot traverse this road, it is forced to produce an unrealistic trajectory containing a u-turn that did not happen, and completely misses the portion of the trajectory on the western north-south road.

Figure~\ref{fig:res3} shows a more extreme failure of the road-constrained map matcher due to a missing road.  In this case, the missing road segment was closed for construction according to the version of the map used by the map matcher, however the road had reopened.  The inability to connect the two halves of the GPS trace here causes the road-constrained map matcher to generate a very implausible trajectory in an attempt to best fit the input data.  The on-/off-road map matcher, in contrast, successfully reconstructs a trajectory very close to that actually taken.


\subsection{Map Learning}
\label{sec:results-map-learning}

The results in Figs.~ \ref{fig:res1},~\ref{fig:res1b}, ~\ref{fig:res2} and ~\ref{fig:res3} show how the sIMM map matcher can produce trajectories that are robust to the class of map errors that would make an on-road-only map matcher generate aberrant trajectories: missing roads, missing connectivity, incorrectly mapped turn restrictions (e.g., a turn restriction that does not exist in the physical world, but is mapped in our digital map used for map matching) or incorrectly mapped road directions (e.g., a road direction that is mapped in the wrong direction in our digital map).
The off-road tracker in the sIMM map matcher is triggered when the map is wrong. Therefore, its output can be used to detect map errors in the digital map in use, and, with sufficient data, to propose corrections. Map errors are ubiquitous in any digital map as a digital map must be constantly changing to match the physical world: new roads, temporary closed or opened roads, permanently closed or opened road, ill-mapped turn restrictions and road directions, vandalized digital map, etc...

\begin{figure}
\includegraphics[height=2in]{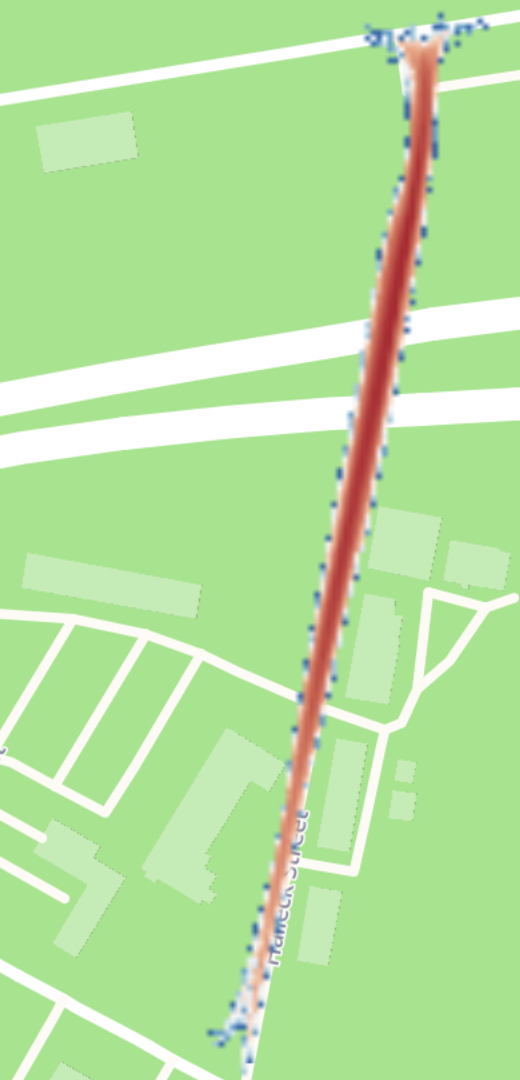}
\includegraphics[height=2in]{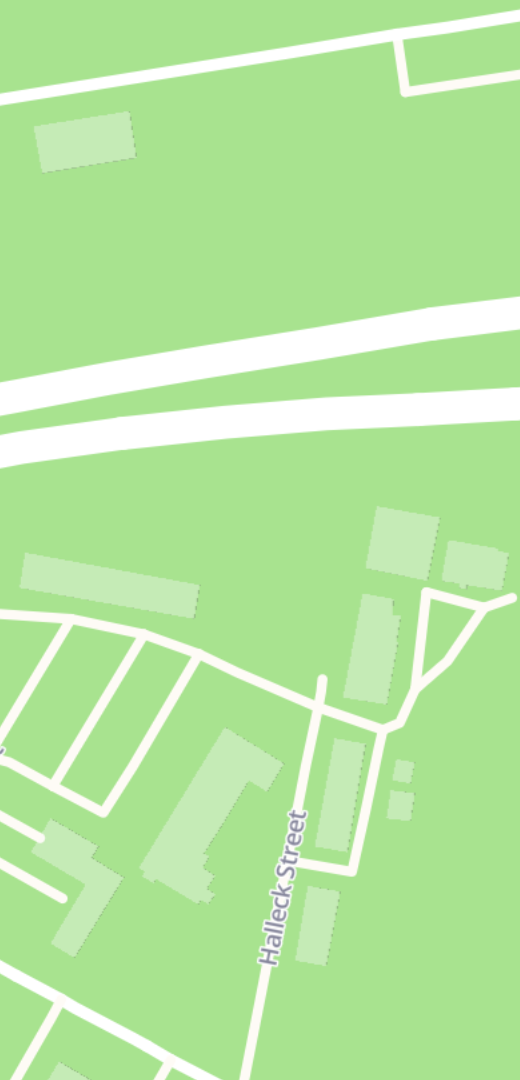}
\includegraphics[height=2in]{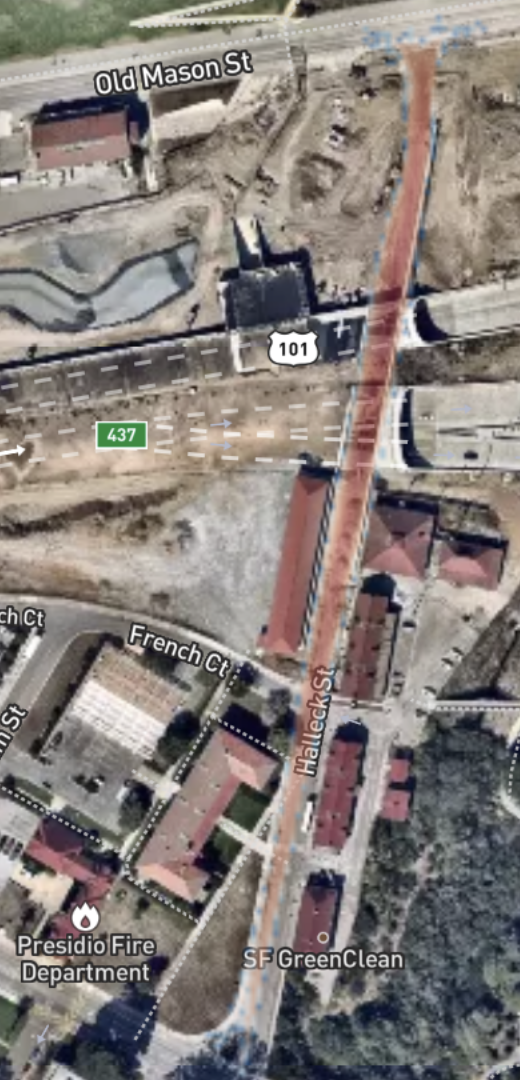}
\caption{After running the sIMM map matcher for one month of Lyft driver location data in early 2019, and using the January 1st, 2018 version of OSM, we detected a large number of off-road traces in the Presidio of San Francisco, displayed by the left panel on which the off-road traces are showed in blue (low density) to red (high density) overlayed on top of OSM. The off-road traces are removed in the middle panel and indeed shows an absence of road. The satellite imagery on the right panel, on which is overlayed in the off-road traces with lower opacity, showed that there were constructions and that this road now exists.}
\label{fig:presidio}
\end{figure}

In the case of passenger vehicle tracking, vehicles generally (but not always!) drive legally and in `drivable' areas (roadways, parking lots, etc.), any recourse to an off-road trajectory could be indicative of a position at which the map used in map matching does not correctly capture true vehicle motion.  Although some apparent off-road motion will be due to illegal manoeuvres (still possibly of interest), poor quality input data (e.g. bad GPS signal) or temporary features (for example, during construction), given a sufficiently large collection of driver traces, map-errors will become apparent as consistent, repeated portions of off-road trajectories.  If a large number of similar `off-road' trajectory portions are available in a given section of the map, it signals a map error with high confidence. As a consequence, it is possible to propose corrections to the digital map in use by considering the collection of trajectories and the constraints of standard road geometry as shown in Fig.~\ref{fig:presidio}, and to assess map quality using large-scale GPS traces.

\subsection{Performance}
\label{sec:results-performance}
The runtime of the sIMM trajectory generation (filter, followed by backward sampling) is around two to four times that of the the standard HMM road-constrained map matcher.  This is due to the additional complexity of the filter and backward sampling, and the fact that some optimizations around pruning of zero-probability branches in the HMM implementation are no longer possible when off-road motion is allowed to be considered.  The implementation of the sIMM filter and backward sampling used here is also not as optimized as that of the standard road-constrained method, so further runtime improvement is probably possible.

\section{Conclusions}
\label{sec:conclusions}

This work has shown how sample-based (e.g. HMM) map matching systems can be made robust to map errors in an efficient and principled way, by combining them with a closed-form free space filter.  In contrast to existing interacting multiple model approaches, the approximation of `semi-interaction' has been adopted in order to allow computation efficient enough for web-scale applications.  This assumption allows a single free-space filter to run without requiring input from the on-road tracking filter, but allows the on-road tracking system to fall back to the free-space (off-road) filter output when needed.  This has applications in improving the robustness of map matching systems by avoiding reconstruction errors when maps contain missing or incorrect data.  Furthermore, with sufficient data, such a system can be used to discover and even propose corrections to errors in the underlying map.

As part of the system a forward filter was developed that is suitable for efficient realtime tracking of vehicles moving primarily on-road but also sometimes off-road or in unmapped areas.  

The method introduced here is general, in the sense that it does not rely on specific on-road or off-road tracking models or filter designs, allowing a range of sample-based on-road trackers and closed-form (or approximately so) off-road trackers to be used.

\bibliographystyle{ACM-Reference-Format}
\bibliography{refs}


\begin{thebibliography}{31}


\ifx \showCODEN    \undefined \def \showCODEN     #1{\unskip}     \fi
\ifx \showDOI      \undefined \def \showDOI       #1{#1}\fi
\ifx \showISBNx    \undefined \def \showISBNx     #1{\unskip}     \fi
\ifx \showISBNxiii \undefined \def \showISBNxiii  #1{\unskip}     \fi
\ifx \showISSN     \undefined \def \showISSN      #1{\unskip}     \fi
\ifx \showLCCN     \undefined \def \showLCCN      #1{\unskip}     \fi
\ifx \shownote     \undefined \def \shownote      #1{#1}          \fi
\ifx \showarticletitle \undefined \def \showarticletitle #1{#1}   \fi
\ifx \showURL      \undefined \def \showURL       {\relax}        \fi
\providecommand\bibfield[2]{#2}
\providecommand\bibinfo[2]{#2}
\providecommand\natexlab[1]{#1}
\providecommand\showeprint[2][]{arXiv:#2}

\bibitem[\protect\citeauthoryear{Arulampalam, Maskell, Gordon, and
  Clapp}{Arulampalam et~al\mbox{.}}{2002}]%
        {Arulampalam2002}
\bibfield{author}{\bibinfo{person}{M~Sanjeev Arulampalam},
  \bibinfo{person}{Simon Maskell}, \bibinfo{person}{Neil Gordon}, {and}
  \bibinfo{person}{Tim Clapp}.} \bibinfo{year}{2002}\natexlab{}.
\newblock \showarticletitle{A tutorial on particle filters for online
  nonlinear/non-Gaussian Bayesian tracking}.
\newblock \bibinfo{journal}{\emph{IEEE Transactions on signal processing}}
  \bibinfo{volume}{50}, \bibinfo{number}{2} (\bibinfo{year}{2002}),
  \bibinfo{pages}{174--188}.
\newblock


\bibitem[\protect\citeauthoryear{Bastani, He, Alizadeh, Balakrishnan, Madden,
  Chawla, Abbar, and DeWitt}{Bastani et~al\mbox{.}}{2018}]%
        {Bastani.2018.cvpr}
\bibfield{author}{\bibinfo{person}{Favyen Bastani}, \bibinfo{person}{Songtao
  He}, \bibinfo{person}{Mohammad Alizadeh}, \bibinfo{person}{Hari
  Balakrishnan}, \bibinfo{person}{Samuel Madden}, \bibinfo{person}{Sanjay
  Chawla}, \bibinfo{person}{Sofiane Abbar}, {and} \bibinfo{person}{David
  DeWitt}.} \bibinfo{year}{2018}\natexlab{}.
\newblock \showarticletitle{{RoadTracer: Automatic Extraction of Road Networks
  from Aerial Images}}. In \bibinfo{booktitle}{\emph{Computer Vision and
  Pattern Recognition (CVPR)}}. \bibinfo{address}{Salt Lake City, UT}.
\newblock


\bibitem[\protect\citeauthoryear{Biagioni and Eriksson}{Biagioni and
  Eriksson}{2012a}]%
        {Biagioni.2012.trr}
\bibfield{author}{\bibinfo{person}{James Biagioni} {and} \bibinfo{person}{Jakob
  Eriksson}.} \bibinfo{year}{2012}\natexlab{a}.
\newblock \showarticletitle{Inferring Road Maps from Global Positioning System
  Traces: Survey and Comparative Evaluation}.
\newblock \bibinfo{journal}{\emph{Transportation Research Record}}
  \bibinfo{volume}{2291}, \bibinfo{number}{1} (\bibinfo{year}{2012}),
  \bibinfo{pages}{61--71}.
\newblock
\urldef\tempurl%
\url{https://doi.org/10.3141/2291-08}
\showDOI{\tempurl}


\bibitem[\protect\citeauthoryear{Biagioni and Eriksson}{Biagioni and
  Eriksson}{2012b}]%
        {Biagioni.2012.sigspatial}
\bibfield{author}{\bibinfo{person}{James Biagioni} {and} \bibinfo{person}{Jakob
  Eriksson}.} \bibinfo{year}{2012}\natexlab{b}.
\newblock \showarticletitle{Map Inference in the Face of Noise and Disparity}.
  In \bibinfo{booktitle}{\emph{Proceedings of the 20th International Conference
  on Advances in Geographic Information Systems}}
  \emph{(\bibinfo{series}{SIGSPATIAL '12})}. \bibinfo{publisher}{ACM},
  \bibinfo{address}{New York, NY, USA}, \bibinfo{pages}{79--88}.
\newblock
\showISBNx{978-1-4503-1691-0}
\urldef\tempurl%
\url{https://doi.org/10.1145/2424321.2424333}
\showDOI{\tempurl}


\bibitem[\protect\citeauthoryear{Blom and Bar-Shalom}{Blom and
  Bar-Shalom}{1988}]%
        {Blom1988}
\bibfield{author}{\bibinfo{person}{Henk~AP Blom} {and} \bibinfo{person}{Yaakov
  Bar-Shalom}.} \bibinfo{year}{1988}\natexlab{}.
\newblock \showarticletitle{The interacting multiple model algorithm for
  systems with Markovian switching coefficients}.
\newblock \bibinfo{journal}{\emph{IEEE transactions on Automatic Control}}
  \bibinfo{volume}{33}, \bibinfo{number}{8} (\bibinfo{year}{1988}),
  \bibinfo{pages}{780--783}.
\newblock


\bibitem[\protect\citeauthoryear{Brakatsoulas, Pfoser, Salas, and
  Wenk}{Brakatsoulas et~al\mbox{.}}{2005}]%
        {Brakatsoulas2005}
\bibfield{author}{\bibinfo{person}{Sotiris Brakatsoulas},
  \bibinfo{person}{Dieter Pfoser}, \bibinfo{person}{Randall Salas}, {and}
  \bibinfo{person}{Carola Wenk}.} \bibinfo{year}{2005}\natexlab{}.
\newblock \showarticletitle{On map-matching vehicle tracking data}. In
  \bibinfo{booktitle}{\emph{Proceedings of the 31st international conference on
  Very large data bases}}. VLDB Endowment, \bibinfo{pages}{853--864}.
\newblock


\bibitem[\protect\citeauthoryear{Capp{\'e}, Godsill, and Moulines}{Capp{\'e}
  et~al\mbox{.}}{2007}]%
        {Cappe2007}
\bibfield{author}{\bibinfo{person}{Olivier Capp{\'e}}, \bibinfo{person}{Simon~J
  Godsill}, {and} \bibinfo{person}{Eric Moulines}.}
  \bibinfo{year}{2007}\natexlab{}.
\newblock \showarticletitle{An overview of existing methods and recent advances
  in sequential Monte Carlo}.
\newblock \bibinfo{journal}{\emph{Proc. IEEE}} \bibinfo{volume}{95},
  \bibinfo{number}{5} (\bibinfo{year}{2007}), \bibinfo{pages}{899--924}.
\newblock


\bibitem[\protect\citeauthoryear{Chen, Lu, Huang, Yang, Gunopulos, and
  Guibas}{Chen et~al\mbox{.}}{2016}]%
        {Chen.2016.kdd}
\bibfield{author}{\bibinfo{person}{Chen Chen}, \bibinfo{person}{Cewu Lu},
  \bibinfo{person}{Qixing Huang}, \bibinfo{person}{Qiang Yang},
  \bibinfo{person}{Dimitrios Gunopulos}, {and} \bibinfo{person}{Leonidas
  Guibas}.} \bibinfo{year}{2016}\natexlab{}.
\newblock \showarticletitle{City-Scale Map Creation and Updating Using GPS
  Collections}. In \bibinfo{booktitle}{\emph{Proceedings of the 22Nd ACM SIGKDD
  International Conference on Knowledge Discovery and Data Mining}}
  \emph{(\bibinfo{series}{KDD '16})}. \bibinfo{publisher}{ACM},
  \bibinfo{address}{New York, NY, USA}, \bibinfo{pages}{1465--1474}.
\newblock
\showISBNx{978-1-4503-4232-2}
\urldef\tempurl%
\url{https://doi.org/10.1145/2939672.2939833}
\showDOI{\tempurl}


\bibitem[\protect\citeauthoryear{Chen and Bierlaire}{Chen and
  Bierlaire}{2015}]%
        {Chen2015}
\bibfield{author}{\bibinfo{person}{Jingmin Chen} {and} \bibinfo{person}{Michel
  Bierlaire}.} \bibinfo{year}{2015}\natexlab{}.
\newblock \showarticletitle{Probabilistic multimodal map matching with rich
  smartphone data}.
\newblock \bibinfo{journal}{\emph{Journal of Intelligent Transportation
  Systems}} \bibinfo{volume}{19}, \bibinfo{number}{2} (\bibinfo{year}{2015}),
  \bibinfo{pages}{134--148}.
\newblock


\bibitem[\protect\citeauthoryear{Cheng and Singh}{Cheng and Singh}{2007}]%
        {Cheng2007}
\bibfield{author}{\bibinfo{person}{Yang Cheng} {and} \bibinfo{person}{Tarunraj
  Singh}.} \bibinfo{year}{2007}\natexlab{}.
\newblock \showarticletitle{Efficient particle filtering for road-constrained
  target tracking}.
\newblock \bibinfo{journal}{\emph{IEEE Trans. Aerospace Electron. Systems}}
  \bibinfo{volume}{43}, \bibinfo{number}{4} (\bibinfo{year}{2007}).
\newblock


\bibitem[\protect\citeauthoryear{Doucet, De~Freitas, and Gordon}{Doucet
  et~al\mbox{.}}{2001}]%
        {Doucet2001}
\bibfield{author}{\bibinfo{person}{Arnaud Doucet}, \bibinfo{person}{Nando
  De~Freitas}, {and} \bibinfo{person}{Neil Gordon}.}
  \bibinfo{year}{2001}\natexlab{}.
\newblock \showarticletitle{An introduction to sequential Monte Carlo methods}.
\newblock In \bibinfo{booktitle}{\emph{Sequential Monte Carlo methods in
  practice}}. \bibinfo{publisher}{Springer}, \bibinfo{pages}{3--14}.
\newblock


\bibitem[\protect\citeauthoryear{Goh, Dauwels, Mitrovic, Asif, Oran, and
  Jaillet}{Goh et~al\mbox{.}}{2012}]%
        {Goh2012}
\bibfield{author}{\bibinfo{person}{Chong~Yang Goh}, \bibinfo{person}{Justin
  Dauwels}, \bibinfo{person}{Nikola Mitrovic}, \bibinfo{person}{Muhammad~Tayyab
  Asif}, \bibinfo{person}{Ali Oran}, {and} \bibinfo{person}{Patrick Jaillet}.}
  \bibinfo{year}{2012}\natexlab{}.
\newblock \showarticletitle{Online map-matching based on hidden markov model
  for real-time traffic sensing applications}. In
  \bibinfo{booktitle}{\emph{Intelligent Transportation Systems (ITSC), 2012
  15th International IEEE Conference on}}. IEEE, \bibinfo{pages}{776--781}.
\newblock


\bibitem[\protect\citeauthoryear{Goyal and Yuen}{Goyal and Yuen}{2019}]%
        {Goyal.2019.kdd}
\bibfield{author}{\bibinfo{person}{Deeksha Goyal} {and} \bibinfo{person}{Albert
  Yuen}.} \bibinfo{year}{2019}\natexlab{}.
\newblock \showarticletitle{Traffic Control Elements Inference using Telemetry
  Data and Convolutional Neural Networks}. In \bibinfo{booktitle}{\emph{The 8th
  International Workshop on Urban Computing (UrbComp 2019), SIGKDD 2019
  Workshop, to be published}}.
\newblock


\bibitem[\protect\citeauthoryear{He, Bastani, Abbar, Alizadeh, Balakrishnan,
  Chawla, and Madden}{He et~al\mbox{.}}{2018}]%
        {He.2018.sigspatial}
\bibfield{author}{\bibinfo{person}{Songtao He}, \bibinfo{person}{Favyen
  Bastani}, \bibinfo{person}{Sofiane Abbar}, \bibinfo{person}{Mohammad
  Alizadeh}, \bibinfo{person}{Hari Balakrishnan}, \bibinfo{person}{Sanjay
  Chawla}, {and} \bibinfo{person}{Sam Madden}.}
  \bibinfo{year}{2018}\natexlab{}.
\newblock \showarticletitle{RoadRunner: Improving the Precision of Road Network
  Inference from GPS Trajectories}. In \bibinfo{booktitle}{\emph{Proceedings of
  the 26th ACM SIGSPATIAL International Conference on Advances in Geographic
  Information Systems}} \emph{(\bibinfo{series}{SIGSPATIAL '18})}.
  \bibinfo{publisher}{ACM}, \bibinfo{address}{New York, NY, USA},
  \bibinfo{pages}{3--12}.
\newblock
\showISBNx{978-1-4503-5889-7}
\urldef\tempurl%
\url{https://doi.org/10.1145/3274895.3274974}
\showDOI{\tempurl}


\bibitem[\protect\citeauthoryear{Jagadeesh and Srikanthan}{Jagadeesh and
  Srikanthan}{2017}]%
        {Jagadeesh2017}
\bibfield{author}{\bibinfo{person}{George~R Jagadeesh} {and}
  \bibinfo{person}{Thambipillai Srikanthan}.} \bibinfo{year}{2017}\natexlab{}.
\newblock \showarticletitle{Online map-matching of noisy and sparse location
  data with hidden markov and route choice models}.
\newblock \bibinfo{journal}{\emph{IEEE Transactions on Intelligent
  Transportation Systems}} \bibinfo{volume}{18}, \bibinfo{number}{9}
  (\bibinfo{year}{2017}), \bibinfo{pages}{2423--2434}.
\newblock


\bibitem[\protect\citeauthoryear{Julier and Uhlmann}{Julier and
  Uhlmann}{1997}]%
        {Julier1997}
\bibfield{author}{\bibinfo{person}{Simon~J Julier} {and}
  \bibinfo{person}{Jeffrey~K Uhlmann}.} \bibinfo{year}{1997}\natexlab{}.
\newblock \showarticletitle{New extension of the Kalman filter to nonlinear
  systems}. In \bibinfo{booktitle}{\emph{Signal processing, sensor fusion, and
  target recognition VI}}, Vol.~\bibinfo{volume}{3068}. International Society
  for Optics and Photonics, \bibinfo{pages}{182--194}.
\newblock


\bibitem[\protect\citeauthoryear{Luxen and Vetter}{Luxen and Vetter}{2011}]%
        {OSRM}
\bibfield{author}{\bibinfo{person}{Dennis Luxen} {and}
  \bibinfo{person}{Christian Vetter}.} \bibinfo{year}{2011}\natexlab{}.
\newblock \showarticletitle{Real-time routing with OpenStreetMap data}. In
  \bibinfo{booktitle}{\emph{Proceedings of the 19th ACM SIGSPATIAL
  International Conference on Advances in Geographic Information Systems}}
  \emph{(\bibinfo{series}{GIS '11})}. \bibinfo{publisher}{ACM},
  \bibinfo{address}{New York, NY, USA}, \bibinfo{pages}{513--516}.
\newblock
\showISBNx{978-1-4503-1031-4}
\urldef\tempurl%
\url{https://doi.org/10.1145/2093973.2094062}
\showDOI{\tempurl}


\bibitem[\protect\citeauthoryear{Mattyus, Luo, and Urtasun}{Mattyus
  et~al\mbox{.}}{2017}]%
        {Mattyus.2017.iccv}
\bibfield{author}{\bibinfo{person}{Gellert Mattyus}, \bibinfo{person}{Wenjie
  Luo}, {and} \bibinfo{person}{Raquel Urtasun}.}
  \bibinfo{year}{2017}\natexlab{}.
\newblock \showarticletitle{DeepRoadMapper: Extracting Road Topology From
  Aerial Images}. In \bibinfo{booktitle}{\emph{The IEEE International
  Conference on Computer Vision (ICCV)}}.
\newblock


\bibitem[\protect\citeauthoryear{Murphy and Godsill}{Murphy and
  Godsill}{2014}]%
        {Murphy2014}
\bibfield{author}{\bibinfo{person}{James Murphy} {and} \bibinfo{person}{Simon
  Godsill}.} \bibinfo{year}{2014}\natexlab{}.
\newblock \showarticletitle{Road-assisted multiple target tracking in clutter}.
  In \bibinfo{booktitle}{\emph{Information Fusion (FUSION), 2014 17th
  International Conference on}}. IEEE, \bibinfo{pages}{1--8}.
\newblock


\bibitem[\protect\citeauthoryear{Newson and Krumm}{Newson and Krumm}{2009}]%
        {Newson2009}
\bibfield{author}{\bibinfo{person}{Paul Newson} {and} \bibinfo{person}{John
  Krumm}.} \bibinfo{year}{2009}\natexlab{}.
\newblock \showarticletitle{Hidden Markov map matching through noise and
  sparseness}. In \bibinfo{booktitle}{\emph{Proceedings of the 17th ACM
  SIGSPATIAL international conference on advances in geographic information
  systems}}. ACM, \bibinfo{pages}{336--343}.
\newblock


\bibitem[\protect\citeauthoryear{{OpenStreetMap contributors}}{{OpenStreetMap
  contributors}}{2017}]%
        {OpenStreetMap}
\bibfield{author}{\bibinfo{person}{{OpenStreetMap contributors}}.}
  \bibinfo{year}{2017}\natexlab{}.
\newblock \bibinfo{title}{{Planet dump retrieved from https://planet.osm.org
  }}.
\newblock \bibinfo{howpublished}{\url{ https://www.openstreetmap.org }}.
\newblock


\bibitem[\protect\citeauthoryear{Orguner, Schon, and Gustafsson}{Orguner
  et~al\mbox{.}}{2009}]%
        {Orguner2009}
\bibfield{author}{\bibinfo{person}{Umut Orguner}, \bibinfo{person}{Thomas~B
  Schon}, {and} \bibinfo{person}{Fredrik Gustafsson}.}
  \bibinfo{year}{2009}\natexlab{}.
\newblock \showarticletitle{Improved target tracking with road network
  information}. In \bibinfo{booktitle}{\emph{Aerospace conference, 2009 IEEE}}.
  IEEE, \bibinfo{pages}{1--11}.
\newblock


\bibitem[\protect\citeauthoryear{Raymond, Morimura, Osogami, and
  Hirosue}{Raymond et~al\mbox{.}}{2012}]%
        {Raymond2012}
\bibfield{author}{\bibinfo{person}{Rudy Raymond}, \bibinfo{person}{Tetsuro
  Morimura}, \bibinfo{person}{Takayuki Osogami}, {and} \bibinfo{person}{Noriaki
  Hirosue}.} \bibinfo{year}{2012}\natexlab{}.
\newblock \showarticletitle{Map matching with hidden Markov model on sampled
  road network}. In \bibinfo{booktitle}{\emph{Pattern Recognition (ICPR), 2012
  21st International Conference on}}. IEEE, \bibinfo{pages}{2242--2245}.
\newblock


\bibitem[\protect\citeauthoryear{S{\"a}rkk{\"a}}{S{\"a}rkk{\"a}}{2008}]%
        {Sarkka2008}
\bibfield{author}{\bibinfo{person}{Simo S{\"a}rkk{\"a}}.}
  \bibinfo{year}{2008}\natexlab{}.
\newblock \showarticletitle{Unscented Rauch--Tung--Striebel Smoother}.
\newblock \bibinfo{journal}{\emph{IEEE Trans. Automat. Control}}
  \bibinfo{volume}{53}, \bibinfo{number}{3} (\bibinfo{year}{2008}),
  \bibinfo{pages}{845--849}.
\newblock


\bibitem[\protect\citeauthoryear{S{\"a}rkk{\"a}, Vehtari, and
  Lampinen}{S{\"a}rkk{\"a} et~al\mbox{.}}{2007}]%
        {Sarkka2007}
\bibfield{author}{\bibinfo{person}{Simo S{\"a}rkk{\"a}}, \bibinfo{person}{Aki
  Vehtari}, {and} \bibinfo{person}{Jouko Lampinen}.}
  \bibinfo{year}{2007}\natexlab{}.
\newblock \showarticletitle{Rao-Blackwellized particle filter for multiple
  target tracking}.
\newblock \bibinfo{journal}{\emph{Information Fusion}} \bibinfo{volume}{8},
  \bibinfo{number}{1} (\bibinfo{year}{2007}), \bibinfo{pages}{2--15}.
\newblock


\bibitem[\protect\citeauthoryear{Stanojevic, Abbar, Thirumuruganathan, Chawla,
  Filali, and Aleimat}{Stanojevic et~al\mbox{.}}{2017}]%
        {Stanojevic.2017.arvix}
\bibfield{author}{\bibinfo{person}{Rade Stanojevic}, \bibinfo{person}{Sofiane
  Abbar}, \bibinfo{person}{Saravanan Thirumuruganathan},
  \bibinfo{person}{Sanjay Chawla}, \bibinfo{person}{Fethi Filali}, {and}
  \bibinfo{person}{Ahid Aleimat}.} \bibinfo{year}{2017}\natexlab{}.
\newblock \showarticletitle{Kharita: Robust Map Inference using Graph
  Spanners}.
\newblock \bibinfo{journal}{\emph{CoRR}}  \bibinfo{volume}{abs/1702.06025}
  (\bibinfo{year}{2017}).
\newblock
\showeprint[arxiv]{1702.06025}
\urldef\tempurl%
\url{http://arxiv.org/abs/1702.06025}
\showURL{%
\tempurl}


\bibitem[\protect\citeauthoryear{Thiagarajan, Ravindranath, LaCurts, Madden,
  Balakrishnan, Toledo, and Eriksson}{Thiagarajan et~al\mbox{.}}{2009}]%
        {Thiagarajan2009}
\bibfield{author}{\bibinfo{person}{Arvind Thiagarajan}, \bibinfo{person}{Lenin
  Ravindranath}, \bibinfo{person}{Katrina LaCurts}, \bibinfo{person}{Samuel
  Madden}, \bibinfo{person}{Hari Balakrishnan}, \bibinfo{person}{Sivan Toledo},
  {and} \bibinfo{person}{Jakob Eriksson}.} \bibinfo{year}{2009}\natexlab{}.
\newblock \showarticletitle{VTrack: accurate, energy-aware road traffic delay
  estimation using mobile phones}. In \bibinfo{booktitle}{\emph{Proceedings of
  the 7th ACM conference on embedded networked sensor systems}}. ACM,
  \bibinfo{pages}{85--98}.
\newblock


\bibitem[\protect\citeauthoryear{Ulmke and Koch}{Ulmke and Koch}{2006}]%
        {Ulmke2006}
\bibfield{author}{\bibinfo{person}{Martin Ulmke} {and}
  \bibinfo{person}{Wolfgang Koch}.} \bibinfo{year}{2006}\natexlab{}.
\newblock \showarticletitle{Road-map assisted ground moving target tracking}.
\newblock \bibinfo{journal}{\emph{IEEE transactions on Aerospace and Electronic
  Systems}} \bibinfo{volume}{42}, \bibinfo{number}{4} (\bibinfo{year}{2006}).
\newblock


\bibitem[\protect\citeauthoryear{Van Der~Merwe, Doucet, De~Freitas, and
  Wan}{Van Der~Merwe et~al\mbox{.}}{2001}]%
        {VanDerMerwe2001}
\bibfield{author}{\bibinfo{person}{Rudolph Van Der~Merwe},
  \bibinfo{person}{Arnaud Doucet}, \bibinfo{person}{Nando De~Freitas}, {and}
  \bibinfo{person}{Eric~A Wan}.} \bibinfo{year}{2001}\natexlab{}.
\newblock \showarticletitle{The unscented particle filter}. In
  \bibinfo{booktitle}{\emph{Advances in neural information processing
  systems}}. \bibinfo{pages}{584--590}.
\newblock


\bibitem[\protect\citeauthoryear{Wan and Van Der~Merwe}{Wan and Van
  Der~Merwe}{2000}]%
        {Wan2000}
\bibfield{author}{\bibinfo{person}{Eric~A Wan} {and} \bibinfo{person}{Rudolph
  Van Der~Merwe}.} \bibinfo{year}{2000}\natexlab{}.
\newblock \showarticletitle{The unscented Kalman filter for nonlinear
  estimation}. In \bibinfo{booktitle}{\emph{Adaptive Systems for Signal
  Processing, Communications, and Control Symposium 2000. AS-SPCC. The IEEE
  2000}}. Ieee, \bibinfo{pages}{153--158}.
\newblock


\bibitem[\protect\citeauthoryear{White, Bernstein, and Kornhauser}{White
  et~al\mbox{.}}{2000}]%
        {White2000}
\bibfield{author}{\bibinfo{person}{Christopher~E White}, \bibinfo{person}{David
  Bernstein}, {and} \bibinfo{person}{Alain~L Kornhauser}.}
  \bibinfo{year}{2000}\natexlab{}.
\newblock \showarticletitle{Some map matching algorithms for personal
  navigation assistants}.
\newblock \bibinfo{journal}{\emph{Transportation research part c: emerging
  technologies}} \bibinfo{volume}{8}, \bibinfo{number}{1-6}
  (\bibinfo{year}{2000}), \bibinfo{pages}{91--108}.
\newblock


\end{thebibliography}

\end{document}